%% file: main.tex
\documentclass[a4paper,12pt]{article}

\usepackage[a4paper,
mag=1000, includefoot,
left=3cm, right=1.5cm, top=2cm, bottom=2cm, headsep=1cm, footskip=1cm]{geometry}
\usepackage[utf8]{inputenc}
\usepackage[english]{babel}

\usepackage{amsmath}
\usepackage{amssymb}
\usepackage{amsfonts}
\usepackage{amsthm}
\usepackage{bbm}
\usepackage{indentfirst}
\usepackage{yhmath}
\usepackage{stackrel}

\newcommand{\R}{\mathbb{R}}
\newcommand{\E}{\mathbf{E}}
\newcommand{\D}{\mathcal{D}}
\newcommand{\A}{\mathcal{A}}
\newcommand{\F}{\mathcal{F}}
\newcommand\Item[1][]{%
  \ifx\relax#1\relax  \item \else \item[#1] \fi
  \abovedisplayskip=0pt\abovedisplayshortskip=0pt~\vspace*{-\baselineskip}}
\newcommand\blfootnote[1]{%
  \begingroup
  \renewcommand\thefootnote{}\footnote{#1}%
  \addtocounter{footnote}{-1}%
  \endgroup
}

\newtheorem{lemma}{Lemma}[section]
\newtheorem{theorem}{Theorem}[section]
\newtheorem*{corollary*}{Corollary}
\newtheorem*{remark*}{Remark}

\title{
    \textbf{Limit theorems for the generator of a symmetric L\'evy process with the delta potential}
    }
\author{
    \textbf{T.~E.~Abildaev} \\
    \normalsize{St.~Petersburg Department of V.~A.~Steklov} \\
    \normalsize{Mathematical Institute of the Russian Academy of Sciences,} \\
    \normalsize{St.~Petersburg State University} \\
}
\date{}

\begin{document}

\maketitle

\begin{abstract}
We consider a one-dimensional symmetric L\'evy process $\xi(t)$, $t \geq 0$, that has local time, which we denote by $L(t,x)$. In the first part, we construct the operator $\A + \mu \, \delta(x - a)$, $\mu > 0$, where $\A$ is the generator of $\xi(t)$, and $\delta(x - a)$ is the Dirac delta function at $a \in \R$. We show that the constructed operator is the generator of $\{ U_t \}_{t \geq 0}$ -- $C_0$-semigroup on $L_2(\R)$, which is given by
\begin{equation*}
    (U_t f)(x)
    = \E f(x - \xi(t)) e^{\mu L(t,x - a)} \text{,} \quad f \in L_2(\R) \cap C_b(\R) \text{,}
\end{equation*}
and prove the Feynman-Kac formula for the delta function-type potentials. We also prove a limit theorem for $U_t f$. In the second part, we construct the measure
\begin{equation*}
    \mathbf{Q}_{T,x}^\mu
    = \frac{e^{\mu L(T,x - a)}}{\E e^{\mu L(T,x - a)}} \, \mathbf{P}_{T,x} \text{,}
\end{equation*}
where $\mathbf{P}_{T,x}$ is the measure of the process $\xi(t)$, $t \leq T$. We show that this measure weakly converges to a Feller process as $T \to \infty$ and prove a limit theorem for the distribution of $\xi(T)$ under $\mathbf{Q}_{T,x}^\mu$.
\end{abstract}


\newpage

\input{1-introduction}
\input{2-V-2-beta_space}
\input{3-psi-lambda_func}
\input{4-operator}
\input{5-semigroups}
\input{6-attracting_distribution}

\newpage

\bibliographystyle{ieeetr}
\bibliography{bibliography}

\end{document}

%% file: 1-introduction.tex
\section{Introduction}

Consider a one-dimensional symmetric L\'evy process $\xi(t)$, $t\geq 0$. It is well known \cite[ch. 2, 2.4]{Apple} that the characteristic function of $\xi(t)$ is given by the L\'evy-Khinchin representation
\begin{equation*}
    \E e^{ip\xi(t)}
    = e^{-t\Psi(p)} \text{,} \quad \Psi(p)
    = \frac{\sigma^2 p^2}{2}
    + \int \limits_{|y|>0} \big( 1 - \text{cos}(py) \big) \, \Pi(dy) \text{,}
\end{equation*}
where $\sigma^2 \geq 0$, and $\Pi$ is the L\'evy measure of the process $\xi(t)$, that is, a symmetric $\sigma$-finite measure that satisfies the condition
\begin{equation*}
    \int \limits_{|y|>0} \text{min}(1, y^2) \, \Pi(dy) < \infty \text{.}
\end{equation*}

Recall that the local time up to time $t$ of the process $\xi(\tau)$ is the density $L(t,\cdot)$ with respect to the Lebesgue measure of its occupation measure, that is, a random measure $\mu_t \text{:} \ \mu_t(\Gamma)=\text{mes} \, \{ \tau <t \, | \, \xi(\tau) \in \Gamma \}$, $\Gamma\in\mathcal{B}(\R)$, if only this density exists.

In this paper, we assume that the local time of $\xi(t)$ exists. This is equivalent to (\cite[ch. V, 1]{Bertoin}, \cite[ch. I, 4.30]{RogWil}) the condition 
\begin{equation}
\label{eq:psi_L1_cond}
    \int \limits_\R \frac{dp}{1 + \Psi(p)} < \infty \text{.}
\end{equation}

One can show (as Borodin--Ibragimov did in \cite[ch. I, \S 4]{BorIbr} for the stable processes) that
\begin{equation*}
    L(t,x)
    = L_2 \, \text{-} \lim_{M \to \infty} \int \limits_{_M}^M e^{-ipx} \Big( \int \limits_0^t e^{ip \xi(\tau)} \, d\tau \Big) \, dp \text{,}
\end{equation*}
and the local time is continuous in both variables with probability 1. From the continuity of the local time it follows that almost surely
\begin{equation}
\label{eq:local_time_lim_repr}
    L(t,x) = \lim_{\varepsilon \to 0+} \frac{1}{2\varepsilon} \int \limits_0^t \mathbbm{1}_{(-\varepsilon, \varepsilon)} (x - \xi(\tau)) \, d\tau \text{.}
\end{equation}

Using this formula and the fact that a sequence
\begin{equation*}
    \frac{1}{2\varepsilon} \mathbbm{1}_{(-\varepsilon, \varepsilon)}(x)
\end{equation*}
converges to the Dirac delta function $\delta(x)$ as $\varepsilon\to 0+$ in the sense of generalized functions, we may formally write
\begin{equation}
\label{eq:local_time_delta_repr}
    L(t,x) = \int \limits_0^t \delta(x - \xi(\tau)) \, d\tau \text{.}
\end{equation}

Recall that the process $\xi(t)$ defines a family of Markov processes $\{\xi_x(t)\}_{x\in\R}$ \cite[chapter I]{Bertoin}, where $\xi_x(t) = x - \xi(t)$. From \eqref{eq:local_time_lim_repr} it follows that $L(t,x)$ is both the local time of $\xi(t)$ at $x$ and the local time of $\xi_x(t)$ at 0.

The family $\{\xi_x(t)\}$ entails a strongly continuous semigroup of operators in the space $L_2(\R)$ \cite[ch. 2, 2.4.2]{LorHirBetz}, acting on the functions from $L_2(\R) \cap C_b(\R)$ by the rule
\begin{equation*}
    \big( T_t f \big)(x)
    = \E f(\xi_x(t)) \text{.}
\end{equation*}
The generator of $\{T_t\}$ is an operator $\A$ that acts on $f \in \D(\A)$ by the rule
\begin{equation}
\label{def:classic_generator}
    \big( \A f \big)(x)
    = \frac{\sigma^2}{2} \, f''(x)
    + \int \limits_{\R \setminus \{ 0 \} } \big( f(x - y) - f(x) + yf'(x) \mathbbm{1}_{[-1,1]}(y) \big) \, \Pi(dy) \text{.}
\end{equation}

The semigroup-theoretic approach reveals connection between the process $\xi(t)$ and a Cauchy problem
\begin{equation}
\label{eq:cauchy_problem}
    \frac{\partial u}{\partial t}(t,x)
    = (\A u)(t,x) + V(x) u(t,x) \text{,}
    \quad u(0,x) = f(x) \text{,}
\end{equation}
where $f \in L_2(\R)$, and $V$ belongs to a ``nice'' enough class of functions (e.g., $C^\infty_0(\R)$). According to the Feynman-Kac formula \cite[p. I, ch. 3]{LorHirBetz}, the unique solution for this problem is given by
\begin{equation}
\label{eq:feynman_kac_formula}
    u(t,x)
    = \E f(\xi_x(t)) \, e^{\int \limits_0^t V(\xi_x(\tau)) \, d\tau} \text{.}
\end{equation}

In turn, under certain conditions the potential $V$ induces a distribution
\begin{equation}
\label{eq:gibbs_meas}
    \frac{e^{\int \limits_0^T V(\omega(\tau)) \, d\tau}}{\E e^{\int \limits_0^T V(\omega(\tau)) \, d\tau}} \, \mathbf{P}_{T,x}(d\omega) \text{,}
\end{equation}
where $\mathbf{P}_{T,x}$ is the distribution of the process $\xi_x(t)$, $t \leq T$, over the sample paths of $\xi_x(t)$ \cite[h. I, ch. 4]{LorHirBetz}. In celebrated \cite{RoyValYor05}, Roynette--Vallois--Yor study this type of measures and call them penalizing. Informally speaking, such measures impose an exponential penalty on the sample paths of the process. 

In problem \eqref{eq:cauchy_problem}, replace the potential $V$ with $\mu \, \delta(x - a)$, $\mu > 0$, $a \in \R$. Formally, using the Feynman-Kac formula, we can see that the unique solution for this problem is a function
\begin{equation*}
    u(t,x)
    = \E f(\xi_x(t)) \, e^{\mu \int \limits_0^t \delta(\xi_x(\tau) - a) \, d\tau} \text{,}
\end{equation*}
or, if we substitute the exponent according to \eqref{eq:local_time_delta_repr},
\begin{equation}
\label{eq:semi_group_repr}
    u(t,x)
    = \E f(\xi_x(t)) \, e^{\mu L(t,x - a)} \text{.}
\end{equation}
Also, we can see that the generator of a corresponding semigroup is $\A + \mu \, \delta(x - a)$.

At the same time a distribution over the sample paths of $\xi(t)$, $t\leq T$, corresponding to the potential $\mu \, \delta(x - a)$, is given by
\begin{equation}
\label{eq:penal_meas}
    \mathbf{Q}_{T,x}^\mu(d\omega)
    = \frac{e^{\mu L(T,x - a)}}{\E e^{\mu L(T,x - a)}} \, \mathbf{P}_{T,x}(d\omega) \text{.}
\end{equation}
One can say that this distribution penalizes the sample paths of $\xi(t)$ for not visiting the point $a$. In other words, it attracts the sample paths of $\xi(t)$ to $a$.

In this paper, we construct $\A + \mu \, \delta(x - a)$ as a self-adjoint extension of the operator $\A$ so that the extension is the generator of the semigroup of operators corresponding to \eqref{eq:semi_group_repr}. Using the constructed operator, we extend the Feynman-Kac formula to the case of delta function-type potentials and prove a limit theorem for an operator semigroup corresponding to this formula. Furthermore, we construct a one-parameter family of distributions $\{ \mathbf{Q}_{T,x}^\mu \}$ that attract the sample paths of $\xi(t)$ to $a$. We prove a limit theorem for the distribution of a random variable $\xi_x(T)$ with respect to $\mathbf{Q}_{T,x}^\mu$ and show that the distributions $\{ \mathbf{Q}_{T,x}^\mu \}$ weakly converge to the distribution of a Feller process as $T \to \infty$.

Previously, in \cite{CranMolSquar} Cranston--Molchanov--Squartini constructed a generator with the delta potential and the corresponding distribution over the sample paths for the stable processes. In \cite{IbrSmoFad24-1}, \cite{IbrSmoFad24-2} Ibragimov--Smorodina--Faddeev constructed a functional
\begin{equation*}
    \int \limits_0^t q(x - w(\tau)) \, d\tau \text{,}
\end{equation*}
where $w$ is the standard Wiener process, and $q$ is a generalized function satisfying a certain condition, and extended the Feynman-Kac formula to the case of $q$-potential.

The author expresses his deep gratitude to N.~V.~Smorodina for statement of the problem and useful discussions.

%% file: 2-V-2-beta_space.tex
\section{Function space $V_2^\beta(\R)$}

By the Fourier transform, $\F$, of a function $g \in L_2(\R)$ we consider a function $\widehat{g} \in L_2(\R)$, which is given by
\begin{equation}
\label{eq:Fourier_transform}
    \widehat{g}(p)
    = L_2\,\text{-}\lim_{M \to \infty} \int \limits_{_M}^M e^{ipx} g(x) \, dx \text{.}
\end{equation}
The inverse Fourier transform, $\F^{-1}$, of a function $\widehat{g} \in L_2(\R)$ is given by
\begin{equation}
\label{eq:inverse_Fourier_transform}
    \frac{1}{2 \pi} \, L_2 \, \text{-} \lim_{M \to \infty} \int \limits_{_M}^M e^{-ipx} \, \widehat{g}(p) \, dp \text{.}
\end{equation}
Recall that, according to the Carleson's theorem \cite{Lacey}, the limits in \eqref{eq:Fourier_transform} and \eqref{eq:inverse_Fourier_transform} may be considered as pointwise.

We begin the construction of the operator $\A + \mu \, \delta(x - a)$ by defining a function space on which the operator $\A$ acts naturally.

Let $\beta \geq 1/2$. Denote by $V_2^\beta(\R)$ a space of functions from $L_2(\R)$ on which a functional $|\cdot|_\beta$ is finite, where
\begin{equation*}
    |\varphi|^2_\beta
    = \int \limits_\R \big( 1 + \Psi^{2 \beta}(p) \big) |\widehat{\varphi}(p)|^2 \, dp \text{.}
\end{equation*}
One can show that $|\cdot|_\beta$ satisfies all the properties of norm, thus $\big( V_2^\beta(\R) \text{,} \, |\cdot|_\beta \big)$, or simply $V_2^\beta(\R)$, is a normed space.

\begin{theorem}
\label{th21}
    Functions from $V_2^\beta(\R)$ are uniformly continuous, bounded, and vanish at infinity. 
\end{theorem}
\begin{proof}
    Let us first prove that
    \begin{equation*}
        \int \limits_\R \frac{dp}{1 + \Psi^{2 \beta}(p)}  < \infty \text{.}
    \end{equation*}

    We have
    \begin{equation*}
        \int \limits_\R \frac{dp}{1 + \Psi^{2 \beta}(p)}  = \int \limits_{\{ \Psi(p) \leq 1 \}} \frac{dp}{1 + \Psi^{2 \beta}(p)}
        + \int \limits_{\{ \Psi(p) > 1 \}} \frac{dp}{1 + \Psi^{2 \beta}(p)}
    \end{equation*}
    \begin{equation*}
        \leq \text{mes} \, \{ \Psi(p) \leq 1 \}
        + \int \limits_{\{ \Psi(p) > 1 \}} \frac{dp}{1 + \Psi(p)}
    \end{equation*}
    \begin{equation*}
        \leq \, 2 \int \limits_{\{ \Psi(p) \leq 1 \}} \frac{dp}{1 + \Psi(p)}
        + \int \limits_{\{ \Psi(p) > 1 \}} \frac{dp}{1 + \Psi(p)}
        \leq 3 \int \limits_\R \frac{dp}{1 + \Psi(p)} < \infty \text{.}
    \end{equation*}
    
    Now let $\varphi\in V_2^\beta(\R)$. Estimating $L_1$-norm of its Fourier transform, we get
    \begin{equation*}
        \| \widehat{\varphi} \|_1
        = \int \limits_\R |\widehat{\varphi}(p)| \, dp
        = \int \limits_\R \frac{\sqrt{1 + \Psi^{2\beta}(p)}}{\sqrt{1 + \Psi^{2\beta}(p)}} \, |\widehat{\varphi}(p)| \, dp
    \end{equation*}
    \begin{equation*}
        \leq \int \limits_\R \frac{dp}{1 + \Psi^{2\beta}(p)} \, \int \limits_\R \big( 1 + \Psi^{2\beta}(p) \big) |\widehat{\varphi}(p)|^2 \, dp
        = \int \limits_\R \frac{dp}{1 + \Psi^{2\beta}(p)} \, | \varphi |_\beta^2
        < \infty \text{.}
    \end{equation*}
    Thus, $\widehat{\varphi} \in L_1(\R)$. 
    
    Continuity, boundedness, and vanishing at infinity follow from the properties of the Fourier transform of functions from $L_1(\R)$ \cite[ch. 10, \S 5]{MakPod} and the Riemann-Lebesgue lemma.

    Let us prove the uniform continuity. Let $x\text{,}\, y \in\R$. We have
    \begin{equation*}
        |\varphi(x) - \varphi(y)|
        = \bigg| \frac{1}{2 \pi} \int \limits_\R (e^{-ipx} - e^{-ipy}) \, \widehat{\varphi}(p) \, dp \bigg|
    \end{equation*}
    \begin{equation*}
        \leq \frac{1}{2 \pi} \int \limits_\R |e^{-ip(x-y)} - 1| |\widehat{\varphi}(p)| \, dp
        = \frac{2}{\pi} \int \limits_\R \Big| \text{sin} \bigg( \frac{p(x-y)}{2} \bigg) \Big| \, |\widehat{\varphi}(p)| \, dp
    \end{equation*}
    \begin{equation*}
        \leq \frac{2}{\pi} \int \limits_{|p| \leq R} \Big| \text{sin} \bigg( \frac{p(x-y)}{2} \bigg) \Big| \, |\widehat{\varphi}(p)| \, dp
        \, + \, \frac{2}{\pi} \int \limits_{|p| > R} \Big| \text{sin} \bigg( \frac{p(x-y)}{2} \bigg) \Big| \, |\widehat{\varphi}(p)| \, dp
    \end{equation*}
    \begin{equation*}
        = \frac{R \, |x-y|}{\pi} \int \limits_{|p| \leq R} |\widehat{\varphi}(p)| \, dp
        \, + \, \frac{2}{\pi} \int \limits_{|p| > R} |\widehat{\varphi}(p)| \, dp \text{.}
    \end{equation*}

    Choosing $R$ so that
    \begin{equation*}
        \frac{2}{\pi} \int \limits_{|p| > R} |\widehat{\varphi}(p)| \, dp < \frac{\varepsilon}{2} \text{,}
    \end{equation*}
    we get that if $|x-y| < (\pi\varepsilon\| \varphi\|_1)/(2R)$, then
    \begin{equation*}
        |\varphi(x) - \varphi(y)| < \varepsilon \text{,}
    \end{equation*}
    which proves the theorem.
\end{proof}

\begin{theorem}
    The space $V_2^\beta(\R)$ is complete and dense in $L_2(\R)$.
\end{theorem}
\begin{proof}
    Let us start with completeness. Let $\{ u_n \}$ be a fundamental sequence from $V_2^\beta(\R)$. Consider the sequence $\{ w_n \}$,
    \begin{equation*}
        \widehat{w}_n(p)
        = \sqrt{1 + \Psi^{2 \beta}(p)} \, \widehat{u}_n(p) \text{.}
    \end{equation*}
    It is fundamental in $L_2(\R)$, because
    \begin{equation*}
        \|w_n - w_m\|_2^2
        = \frac{1}{2\pi} \int \limits_\R \big| \widehat{w}_n(p) - \widehat{w}_m(p) \big|^2 \, dp
    \end{equation*}
    \begin{equation*}
        = \frac{1}{2\pi} \int \limits_\R \big( 1 + \Psi^{2 \beta}(p) \big) \big| \widehat{u}_n(p) - \widehat{u}_m(p) \big|^2 \, dp
        = \frac{1}{2\pi} \, | u_n - u_m |_\beta^2 \stackrel[n \to \infty]{}{\longrightarrow} 0 \text{.}
    \end{equation*}
    
    Denote by $w$ the limit of $\{ w_n \}$ in $L_2(\R)$ and define a function $u$ through its Fourier transform, assuming 
    \begin{equation*}
        \widehat{u}(p)
        = \frac{\widehat{w}(p)}{\sqrt{1 + \Psi^{2 \beta}(p)}} \text{.}
    \end{equation*}
    
    It is clear that $u \in V_2^\beta(\R)$ and
    \begin{equation*}
        | u - u_n |_\beta^2
        = \int \limits_\R \big( 1 + \Psi^{2 \beta}(p) \big) \big| \widehat{u}(p) - \widehat{u}_n(p) \big|^2 \, dp
        = \int \limits_\R \big| \widehat{w}(p) - \widehat{w}_n(p) \big|^2 \, dp \stackrel[n \to \infty]{}{\longrightarrow} 0 \text{.}
    \end{equation*}
    
    Now let us show that $V_2^\beta(\R)$ is dense in $L_2(\R)$. Let $u \in L_2(\R)$. Consider a sequence $\{ u_n \}$, defined by
    \begin{equation*}
        \widehat{u}_n(p)
        =
        \begin{cases}
            \widehat{u}(p)\text{,} \quad &|p| < n \text{,} \\
            \widehat{u}(p) \big/ \sqrt{1 + \Psi^{2 \beta}(p)} \text{,} \quad &|p| \geq n \text{.}
        \end{cases}
    \end{equation*}
    One can easily see that $u_n \in V_2^\beta(\R)$. Moreover,
    \begin{equation*}
        \| u - u_n \|_2^2
        = \frac{1}{2\pi} \int \limits_\R \big| \widehat{u}(p) - \widehat{u}_n(p) \big|^2 \, dp
        = \frac{1}{2\pi} \int \limits_{|p| \geq n} \big| \widehat{u}(p) - \widehat{u}_n(p) \big|^2 \, dp
    \end{equation*}
    \begin{equation*}
        = \frac{1}{2\pi} \int \limits_{|p| \geq n} \Bigg( 1 - \frac{1}{\sqrt{1 + \Psi^{2 \beta}(p)}} \Bigg)^2 \big| \widehat{u}(p) \big|^2 \, dp
        = \frac{1}{2\pi} \int \limits_{|p| \geq n} |\widehat{u}(p)|^2 \, dp \stackrel[n \to \infty]{}{\longrightarrow} 0 \text{,}
    \end{equation*}
    which completes the proof.
\end{proof}

It is evident that if $u \in \D(\A) \cap V_2^\beta(\R)$, then the formula \eqref{def:classic_generator} takes the form
\begin{equation}
\label{def:adv_classic_generator}
    \big( \A u \big)(x)
    = - \frac{1}{2 \pi} \, L_2 \, \text{-} \lim_{M \to \infty} \int \limits_{-M}^M e^{-ipx} \, \Psi(p) \widehat{u}(p) \, dp \text{.}
\end{equation}
Taking this representation as a basis, we consider the generator of the process $\xi(t)$ as an unbounded, densely defined operator $\A : V_2^1(\R) \rightarrow L_2(\R)$, that acts by the formula \eqref{def:adv_classic_generator}.

From \eqref{def:adv_classic_generator} it follows that the Fourier transform diagonalizes $\A$, that is,
\begin{equation*}
    \F \A = \widehat{\A} \, \F \text{,}
\end{equation*}
where $\widehat{\A}$ is the multiplication operator for the function $-\Psi$.

One can also show that the operator $\A$ is self-adjoint, which allows us for each Borel $f$ to define the operator $f(\A)$ \cite[ch. 6, \S 6.1]{BirSol}. In particular, for $u \in V_2^\beta(\R)$
\begin{equation*}
    \big( (-\A)^\beta u \big)(x)
    \stackrel{\text{def}}{=} \frac{1}{2 \pi} \lim_{M \to \infty} \int \limits_{-M}^M e^{-ipx} \, \Psi^\beta(p) \widehat{u}(p) \, dp \text{.}
\end{equation*}

Let us formulate a lemma that connects the action of the operator $(-\A)^\beta$ with the value of a function at a point.
\begin{lemma}
\label{l21}
Let $u \in V_2^\beta(\R)$, $\varkappa > 0$. Then for any $x \in \R$
    \begin{equation*}
        |u(x)|^2
        \leq \Bigg( \frac{1}{2 \pi} \int \limits_\R \frac{dp}{\varkappa + \Psi^{2 \beta}(p)} \Bigg) \left( \| (-\A)^\beta u \|_2^2 + \varkappa \| u \|_2^2 \right) \text{.}
    \end{equation*}
\end{lemma}
\begin{proof}
    Let $x \in \R$. Using the continuity of $u$, the properties of the limit, and the Schwartz inequality, we get
    \begin{equation*}
        |u(x)|^2
        = \left| \frac{1}{2 \pi} \lim_{M \to \infty} \int \limits_{-M}^M e^{-ipx} \, \widehat{u}(p) \, dp \right|^2
    \end{equation*}
    \begin{equation*}
        = \left| \frac{1}{2 \pi} \lim_{M \to \infty} \int \limits_{-M}^M e^{-ipx} \, \widehat{u}(p) \, \frac{\sqrt{\kappa + \Psi^{2 \beta}(p)}}{\sqrt{\kappa + \Psi^{2 \beta}(p)}} \, dp \right|^2
    \end{equation*}
    \begin{equation*}
        = \lim_{M \to \infty} \left| \frac{1}{2 \pi} \int \limits_{-M}^M e^{-ipx} \, \widehat{u}(p) \, \frac{\sqrt{\kappa + \Psi^{2 \beta}(p)}}{\sqrt{\kappa + \Psi^{2 \beta}(p)}} \, dp \right|^2
    \end{equation*}
    \begin{equation*}
        \leq \lim_{M \to \infty} \Bigg[ \bigg( \frac{1}{2 \pi} \bigg)^2 \int \limits_{-M}^M \frac{dp}{\kappa + \Psi^{2 \beta}(p)} \, \int \limits_{-M}^M \big( \kappa + \Psi^{2 \beta}(p) \big) |\widehat{u}(p)|^2 \, dp \Bigg]
    \end{equation*}
    \begin{equation*}
        = \Bigg( \frac{1}{2 \pi} \int \limits_\R \frac{dp}{\kappa + \Psi^{2 \beta}(p)} \Bigg) \, \frac{1}{2 \pi} \lim_{M \to \infty} \int \limits_{-M}^M \big( \kappa + \Psi^{2 \beta}(p) \big) |\widehat{u}(p)|^2 \, dp
    \end{equation*}
    \begin{equation*}
        = \Bigg( \frac{1}{2 \pi} \int \limits_\R \frac{dp}{\kappa + \Psi^{2 \beta}(p)} \Bigg) \left( \| (-\A)^\beta u \|_2^2 + \kappa \| u \|_2^2 \right) \text{.}
    \end{equation*}
\end{proof}

%% file: 3-psi-lambda_func.tex
\section{Family of functions $\{ \psi_\lambda \}$}

Throughout this section, $\lambda \in \mathbb{C} \setminus (-\infty, 0]$, unless otherwise specified. Let us define a function $\psi_\lambda$ through its Fourier transform, assuming
\begin{equation*}
    \widehat{\psi}_\lambda(p)
    = \frac{1}{\Psi(p) + \lambda} \text{.}
\end{equation*}
From the condition \eqref{eq:psi_L1_cond} it follows that $\widehat{\psi}_\lambda \in L_1(\R)$, therefore
\begin{equation*}
    \psi_\lambda(x)
    = \frac{1}{2 \pi} \int \limits_\R \frac{e^{-ipx}}{\Psi(p) + \lambda} \, dp \text{.}
\end{equation*}

The following statements are about the properties of the function $\psi_\lambda$.

\begin{lemma}
    The function $\psi_\lambda$ is uniformly continuous, bounded, and vanishes at infinity.
\end{lemma}

The proof is similar to the proof of the \ref{th21} theorem.

\begin{lemma}
\label{l23}
Let $\nu > 0$. The function $F(\nu) = \psi_\nu(0)$ is positive, continuous, and monotone decreasing. In addition,
    \begin{equation}
    \label{l23_eq1}
       \lim_{\nu \to 0+} F(\nu) = \infty \text{,} \quad \lim_{\nu \to \infty} F(\nu) = 0.
    \end{equation}
\end{lemma}
\begin{proof}
    Let $0 < \nu_1 < \nu_2$. We have
    \begin{equation*}
       F(\nu_2) - F(\nu_1) = \frac{1}{2 \pi} \int \limits_\R \frac{dp}{(\Psi(p) + \nu_1)(\Psi(p) + \nu_2)} \, (\nu_1 - \nu_2) \text{,}
    \end{equation*}
    which proves continuity and monotone decrease.
    
    Positiveness and equations \eqref{l23_eq1} are obvious.
\end{proof}

The following statement connects the function $\psi_\lambda$ to the local time of $\xi(t)$.
\begin{theorem}
\label{th31}
    Let $\text{Re} \, \lambda > 0$. Then
    \begin{equation*}
        \psi_\lambda(x)
        = \E \int \limits_0^\infty e^{-\lambda t} \, L(dt,x) \text{.}
    \end{equation*}
\end{theorem}
\begin{proof}
    Using the properties of $L(t,x)$, we obtain
    \begin{equation*}
        \E \int \limits_0^\infty e^{-\lambda t} \, L(dt,x)
        = \lim_{T \to \infty} \E \int \limits_0^T e^{-\lambda t} \, L(dt,x)
    \end{equation*}
    \begin{equation*}
        = \lim_{T \to \infty} \E \int \limits_0^T e^{-\lambda t} \, d_t \bigg[ \lim_{M \to \infty} \frac{1}{2 \pi} \int \limits_{-M}^M e^{-ipx} \Big( \int \limits_0^t e^{ip\xi(\tau)} \, d\tau \Big) \, dp \bigg]
    \end{equation*}
    \begin{equation*}
        = \lim_{T \to \infty} \E \int \limits_0^T e^{-\lambda t} \, d_t \bigg[ \lim_{M \to \infty} \frac{1}{2 \pi} \int \limits_0^t \bigg( \int \limits_{-M}^M e^{-ipx} e^{ip\xi(\tau)} \, dp \bigg) \, d\tau \bigg]
    \end{equation*}
    \begin{equation*}
        = \lim_{T \to \infty} \E \int \limits_0^T e^{-\lambda t} \, \bigg( \lim_{M \to \infty} \frac{1}{2 \pi} \int \limits_{-M}^M e^{-ipx} e^{ip\xi(t)} \, dp \bigg) \, dt
    \end{equation*}
    \begin{equation*}
        = \lim_{T \to \infty} \int \limits_0^T e^{-\lambda t} \, \bigg( \frac{1}{2 \pi} \lim_{M \to \infty} \int \limits_{-M}^M e^{-ipx} e^{-t \Psi(p)} \, dp \bigg) \, dt
    \end{equation*}
    \begin{equation*}
        = \frac{1}{2 \pi} \lim_{T \to \infty} \lim_{M \to \infty} \int \limits_{-M}^M e^{-ipx} \bigg( \int \limits_0^T e^{-\lambda t} e^{-t \Psi(p)} \, dt \bigg) \, dp
    \end{equation*}
    \begin{equation*}
        = \frac{1}{2 \pi} \lim_{M \to \infty} \lim_{T \to \infty} \int \limits_{-M}^M e^{-ipx} \, \frac{1 - e^{-(\lambda + \Psi(p))T}}{\Psi(p) + \lambda} \, dp
    \end{equation*}
    \begin{equation*}
        = \frac{1}{2 \pi} \lim_{M \to \infty} \int \limits_{-M}^M \frac{e^{-ipx}}{\Psi(p) + \lambda} \, dp
        = \frac{1}{2 \pi} \int \limits_\R \frac{e^{-ipx}}{\Psi(p) + \lambda} \, dp
        = \psi_\lambda(x) \text{.}
    \end{equation*}
\end{proof}

\begin{lemma}
\label{l24}
Let $\nu > 0$. The function $\psi_\nu$ is even and positive. Besides,
\begin{equation*}
    \| \psi_\nu \|_1
    = \frac{1}{\nu} \text{.}
\end{equation*}
\end{lemma}
\begin{proof}
    Evenness of $\psi_\nu$ follows from evenness of $\Psi$. Positivity of $\psi_\nu$ follows from the theorem \ref{th31}.

    Furthermore,
    \begin{equation*}
        \| \psi_\nu \|_1
        = \int \limits_\R \psi_\nu(x) \, dx
        = \lim_{M \to \infty} \int \limits_{_M}^M e^{ipx} \psi_\nu(x) \, dx \, \Big|_{p = 0}
        = \frac{1}{\Psi(p) + \nu} \, \Big|_{p = 0}
        = \frac{1}{\nu} \text{,}
    \end{equation*}
    which completes the proof.
\end{proof}

Let us formulate and prove the statement about the connection of $\psi_\lambda$ and the resolvent of the operator $\A$.
\begin{theorem}
Let $f\in L_2(\R)\cap L_1(\R)$. Then
\begin{equation*}
    \big( (\A - \lambda)^{-1} f \big)(x)
    = -\int \limits_\R \psi_\lambda(x-y) f(y)  \, dy \text{.}
\end{equation*}
\end{theorem}
\begin{proof}
    It is easy to show that the operator $\A$ resolvent acts by the formula
    \begin{equation*}
        \big( (\A - \lambda)^{-1} f \big)(x)
        = - \frac{1}{2 \pi} \int \limits_\R e^{-ipx} \frac{\widehat{f}(p)}{\Psi(p) + \lambda} \, dp \text{,}
    \end{equation*}
    where
    \begin{equation*}
        \widehat{f}(p)
        = \lim_{M \to \infty} \int \limits_{-M}^M e^{ipx} f(x) \, dx \text{.}
    \end{equation*}

    If $f \in L_2(\R) \cap L_1(\R)$, then by the Fubini's theorem
    \begin{equation*}
        \big( (\A - \lambda)^{-1} f \big)(x)
        = - \frac{1}{2 \pi} \int \limits_\R \frac{e^{-ipx}}{\Psi(p) + \lambda} \bigg( \int \limits_\R e^{ipy} f(y) \, dy \bigg) \, dp
    \end{equation*}
    \begin{equation*}
        = - \int \limits_\R \bigg( \frac{1}{2 \pi} \int \limits_\R \frac{e^{-ip(x - y)}}{\Psi(p) + \lambda} \, dp \bigg) f(y) \, dy
        = -\int \limits_\R \psi_\lambda(x-y) f(y)  \, dy \text{.}
    \end{equation*}
\end{proof}

%% file: 4-operator.tex
\section{Operator $\A + \mu \, \delta(x - a)$}

In this section, we define and study the properties of an operator
\begin{equation*}
    \A_\mu = \A + \mu \, \delta(x - a) \text{,}
\end{equation*}
where $\mu > 0$.

By $\mathcal{D}_0$ denote the function space
\begin{equation*}
    \{ \varphi \in V_2^1(\R) : \varphi(a) = 0 \} \text{.}
\end{equation*}

By $\mathcal{D}(\A_\mu)$ denote the domain of $\A_\mu$ and define it to be
\begin{equation*}
     \mathcal{D}_0 \ \oplus \ \text{Lin}\big( \psi_\nu(\cdot-a) \big) \text{,}
\end{equation*}
where the constant $\nu$ is uniquely determined by the equation
\begin{equation*}
    \mu \, \psi_\nu(0) =
    \frac{\mu}{2\pi} \int \limits_\R \frac{dp}{\Psi(p) + \nu} = 1 \text{.}
\end{equation*}

From Lemma \ref{l23} it follows that $\nu$ is determined correctly.

\begin{remark*}
    Let $\xi(t)$ be a symmetric stable process such that
    \begin{equation*}
        \E e^{ip\xi(t)}
        = e^{-B t |p|^\alpha} \text{,} \quad \alpha \in (1,2] \text{,} \quad B > 0 \text{.}
    \end{equation*}
    Then
    \begin{equation*}
        \nu
        = B^\frac{1}{1 - \alpha} \Bigg( \frac{\mu}{\pi} \int \limits_0^\infty \frac{d\theta}{\theta^\alpha + 1} \Bigg)^{\frac{\alpha}{\alpha - 1}} \text{.}
    \end{equation*}
    In particular, if $\xi(t)$ is the standard Wiener process, then $\nu = \mu^2 / 2$.
\end{remark*}

For $u \in \D(\A_\mu)$, $u = \varphi + C \psi_\nu$, by definition, put
\begin{equation*}
    \A_\mu u
    = \A \varphi + \nu \psi_\nu(\cdot - a) \text{.}
\end{equation*}

Our next aim is to show that the operator $\A_\mu$ is self-adjoint. To do this, we formulate a statement about the action of $\A_\mu$ in terms of the Fourier transform, and then prove that the operator $\A_\mu$ is symmetric and closed. 
\begin{lemma}
\label{l32}
Let $u \in \mathcal{D}(\A_\mu)$. Then
    \begin{equation*}
        \big( \widehat{\A_\mu u} \big)(p)
        = -\Psi(p) \widehat{u}(p) + \mu \, u(a) e^{ipa} \text{.}
    \end{equation*}
\end{lemma}
\begin{proof}
    Let $u = \varphi + C \psi_\nu(\cdot - a)$. We have
    \begin{equation*}
        \big( \widehat{\A_\mu u} \big)(p)
        = \big( \widehat{\A_\mu} \widehat{u} \big)(p)
    \end{equation*}
    \begin{equation*}
        = -\Psi(p) \widehat{\varphi}(p) + \nu \, C e^{ipa} \widehat{\psi}_\nu(p)
        = -\Psi(p) \widehat{\varphi}(p) + C \frac{\nu e^{ipa}}{\Psi(p) + \nu}
    \end{equation*}
    \begin{equation*}
        = -\Psi(p) \widehat{\varphi}(p) + C e^{ipa} - C \, \frac{\Psi(p) e^{ipa}}{\Psi(p) + \nu}
        = -\Psi(p) \widehat{u}(p) + C e^{ipa}
    \end{equation*}
    \begin{equation*}
        = -\Psi(p) \widehat{u}(p) + \mu \, C \psi_\nu(0) e^{ipa}
        = -\Psi(p) \widehat{u}(p) + \mu \, u(a) e^{ipa} \text{,}
    \end{equation*}
    which establishes the formula.
\end{proof}

\begin{theorem}
    The operator $\A_\mu$ is symmetric and closed.
\end{theorem}
\begin{proof}
    Let us begin with symmetricity. Let $u, v \in \mathcal{D}(\A_\mu)$. Using the Lemma \ref{l32}, we obtain
    \begin{equation*}
        (\A_\mu u,v)
        = \frac{1}{2\pi} (\widehat{\A}_\mu \widehat{u}, \widehat{v})
    \end{equation*}
    \begin{equation*}
        = - \frac{1}{2\pi} \big( \Psi \widehat{u}, \widehat{v} \big)
        + \frac{1}{2\pi} \mu \, u(a) (e^{i(\cdot)a}, \widehat{v})
        = - \frac{1}{2\pi} \big( \widehat{u}, \Psi \widehat{v} \big)
        + \mu \, u(a) \overline{v(a)}
    \end{equation*}
    \begin{equation*}
        = - \frac{1}{2\pi} \big( \widehat{u}, \Psi \widehat{v} \big)
        + \frac{1}{2 \pi} \mu \, (\widehat{u}, e^{i(\cdot)a}) \overline{v(a)}
        = \frac{1}{2 \pi} (\widehat{u}, \widehat{\A}_\mu \widehat{v})
        = (u, \A_\mu v) \text{,}
    \end{equation*}
    where $(e^{i(\cdot)a}, \widehat{v})$ and $(\widehat{u}, e^{i(\cdot)a})$ are to be considered as
    \begin{equation*}
         \lim_{M \to \infty} \int \limits_{-M}^M e^{ipa} \, \overline{\widehat{v}(p)} \, dp \ \ \text{and} \ \
         \lim_{M \to \infty} \int \limits_{-M}^M e^{-ipa} \, \widehat{u}(p) \, dp
    \end{equation*}
    respectively.

    Let us proceed with closedness. We need to show that $\D(\A_\mu)$ is complete with respect to the norm $\| \cdot \|_\mu$, where
    \begin{equation*}
        \| g \|_\mu^2
        = \| g \|_2^2 + \| \A_\mu g \|_2^2 \text{,} \quad g \in \D(\A_\mu) \text{.}
    \end{equation*}
    
    Let $\{ u_n \}$ be a sequence from $\D(\A_\mu)$, $u_n = \varphi_n + C_n \psi_\nu(\cdot - a)$, fundamental with respect to $\|\cdot\|_\mu$. It means that
    \begin{align*}
        \| (\varphi_n - \varphi_m) + (C_n - C_m) \psi_\nu(\cdot - a) \|_2 &\to 0 \quad \text{and} \\
        \| \A (\varphi_n - \varphi_m) + \nu (C_n - C_m) \psi_\nu(\cdot - a) \|_2 &\to 0
    \end{align*}
    as $n \text{, } m \to \infty$.

    From the properties of the norm it follows that
    \begin{equation*}
        \| \A (\varphi_n - \varphi_m) + \nu (C_n - C_m) \psi_\nu(\cdot - a) \|_2
    \end{equation*}
    \begin{equation*}
        = \| (\A - \nu) (\varphi_n - \varphi_m) + \nu \big( (\varphi_n - \varphi_m) + (C_n - C_m) \psi_\nu(\cdot - a) \big) \|_2
    \end{equation*}
    \begin{equation*}
        \geq \big| \| (\A - \nu) (\varphi_n - \varphi_m) \|_2 - |\nu| \| (\varphi_n - \varphi_m) + (C_n - C_m) \psi_\nu(\cdot - a) \|_2 \big| \text{,}
    \end{equation*}
    which means that if $n \text{, } m \to \infty$, then
    \begin{equation*}
        \| (\A - \nu) (\varphi_n - \varphi_m) \|_2 \rightarrow 0 \text{,}
    \end{equation*}
    which, in turn, means that if $n \text{, } m \to \infty$, then
    \begin{equation*}
        \| \varphi_n - \varphi_m \|_2 \rightarrow 0 \quad \text{and} \quad \| \A (\varphi_n - \varphi_m) \|_2 \rightarrow 0 \text{.}
    \end{equation*}
    
    Therefore, the fundamentality of $\{ u_n \}$ with respect to the norm of $\|\cdot\|_\mu$ is equivalent to what follows:
    \begin{align*}
        | \varphi_n - \varphi_m |_1^2
        = 2 \pi \big( \| \varphi_n - \varphi_m \|_2^2 + \| \A (\varphi_n - \varphi_m) \|_2^2 \big) &\to 0 \quad \text{and} \\
        |C_n - C_m| &\to 0
    \end{align*}
    as $n \text{,} m \to \infty$.
    
    Due to completeness of $(V_2^1(\R), |\cdot|_1)$ there exists a function $\varphi\in V_2^1(\R)$ such that
    \begin{equation*}
        | \varphi - \varphi_n |_1^2 = 2 \pi \big( \| \varphi - \varphi_n \|_2^2 + \| \A (\varphi - \varphi_n) \|_2^2 \big) \to 0 \text{,} \quad n \to \infty \text{.}
    \end{equation*}
    There is also a constant $C \in \R$ such that
    \begin{equation*}
        |C - C_n| \to 0 \text{,} \quad n \to \infty \text{.}
    \end{equation*}
    
    Thus, if $n \to \infty$ then
    \begin{equation*}
        \| u - u_n \|_\mu \to 0 \text{,}
    \end{equation*}
    and for $u$ to belong to $\D(\A_\mu)$, it is needed that the condition $\varphi(a) = 0$ is met.
    
    Applying Lemma \ref{l21} to $\varphi - \varphi_n$ and assuming $\beta = 1$, $\kappa = 1$, we obtain   
    \begin{equation*}
        |\varphi(a)|^2
        = |\varphi(a) - \varphi_n(a)|^2
    \end{equation*}
    \begin{equation*}
        \leq \Bigg( \frac{1}{2 \pi} \int \limits_\R \frac{dp}{1 + \Psi^2(p)} \Bigg) \big( \| \A (\varphi - \varphi_n) \|_2^2 + \| \varphi - \varphi_n \|_2^2 \big)
    \end{equation*}
    \begin{equation*}
        \leq \Bigg( \frac{1}{2 \pi} \int \limits_\R \frac{dp}{1 + \Psi^2(p)} \Bigg) \frac{1}{2 \pi} \| \varphi - \varphi_n \|_1^2 \to 0 \text{,} \quad n \to \infty \text{,}
    \end{equation*}
    and the proof is complete.
\end{proof}

\begin{theorem}
\label{th42}
    The operator $\A_\mu$ is self-adjoint.
\end{theorem}
\begin{proof}
    Since $\A_\mu$ is symmetric and closed, it is sufficient \cite[ch. 4, \S 4.1]{BirSol} to show that
    \begin{equation}
    \label{th42_proof1}
        \text{Ker}(\A_\mu^* \pm i) = \{ 0 \} \text{.}
    \end{equation}
    
    Suppose, contrary to this, that there exists $v \in \D(\A_\mu^*)$ such that $v \not \equiv 0$ and for any $u\in\D(\A_\mu)$
    \begin{equation*}
        (u, (\A_\mu^* \pm i)v) = 0 \text{.}
    \end{equation*}
    Let $u = \varphi + C \psi_\nu(\cdot - a)$. Then
    \begin{equation*}
        (u, (\A_\mu^* \pm i) v)
        = ((\A_\mu \mp i) u, v)
        = \frac{1}{2 \pi} ((\widehat{\A}_\mu \mp i) \widehat{u}, \widehat{v})
    \end{equation*}
    \begin{equation*}
        = \frac{1}{2 \pi} ((-\Psi \mp i) \widehat{\varphi}, \widehat{v}) + \frac{C}{2 \pi} ((\nu \mp i) \widehat{\psi}_\nu \, e^{i(\cdot)a}, \widehat{v})
        = 0 \text{,}
    \end{equation*}
    which is equivalent to
    \begin{equation*}
        ((-\Psi \mp i) \widehat{\varphi}, \widehat{v})
        = \int \limits_\R (-\Psi(p) \mp i) \widehat{\varphi}(p) \overline{\widehat{v}(p)} \, dp = 0 \text{,}
    \end{equation*}
    \begin{equation}
    \label{th42_proof3}
        (\widehat{\psi}_\nu \, e^{i(\cdot)a}, \widehat{v})
        = \int \limits_\R e^{ipa} \frac{\overline{\widehat{v}(p)}}{\Psi(p) + \nu} \, dp = 0 \text{.}
    \end{equation}
    
    From the last equation it follows that
    \begin{equation*}
        \widehat{v}(p)
        = \frac{C}{\Psi(p) \pm i} \text{.}
    \end{equation*}
    for some $C \in \R$. Substituting the last expression for \eqref{th42_proof3}, we get
    \begin{equation*}
        \int \limits_\R \frac{e^{ipa}}{\Psi(p) \pm i} \, dp
        = \int \limits_\R \frac{e^{ipa}}{\Psi(p) + \nu} \, dp \text{,}
    \end{equation*}
    which is obviously not true for whatever sign before $i$. This means that there is no function $v$ with the claimed properties and, therefore, \eqref{th42_proof1} holds.
\end{proof}

In a Hilbert space $\mathcal{H}$, there is a natural bijection between the semi-bounded from below self-adjoint operators and the closed semi-bounded from below quadratic forms. Using this bijection, we show that $\A_\mu$ is to be considered as the operator $\A + \mu \delta(x - a)$.

Let us recall some concepts. A Hermitian form $a$ defined on a dense subspace of a Hilbert space $\D[a] \subset \mathcal{H}$ is called \textit{semi-bounded from below} if for some $m_a \in \R$ and any $u \in \D[a]$
\begin{equation*}
    a[u,u] \geq m_a \| u \|_\mathcal{H}^2 \text{.}
\end{equation*}
A self-adjoint operator $A$ is called \textit{semi-bounded from below} if the form generated by this operator is semi-bounded from below, that is
\begin{equation*}
    (Au,u) \geq m_a \| u \|_\mathcal{H}^2
\end{equation*}
for some $m_a \in \R$.

Without loss of generality, we assume that $m_a < 0$. A semi-bounded from below form $a$ is \textit{closed} if $\D[a]$ is complete with respect to the norm $\|\cdot\|_a$, where
\begin{equation*}
    \| u \|_a^2 =
    a[u,u] + (-m_a + 1) \| u \|_\mathcal{H}^2 \text{.}
\end{equation*}

A self-adjoint operator $A$ and a Hermitian form $a$ \textit{correspond} to each other if
\begin{equation}
\label{eq:op_and_form_cond1}
    \D(A) \subset \D[a]
\end{equation}
and for any $u \text{, } v \in \D[a]$
\begin{equation}
\label{eq:op_and_form_cond2}
    (Au,v) = a[u,v] \text{.}
\end{equation}

The following statement holds in any Hilbert space.
\begin{itemize}
    \item[a)] Each semi-bounded from below self-adjoint operator corresponds to the unique closed semi-bounded from below Hermitian form.
    \item[b)] Each closed semi-bounded from below Hermitian form corresponds to the unique semi-bounded from below self-adjoint operator.
\end{itemize}

Let us now proceed with the construction of the form $a_\mu$ that corresponds to the semi-bounded self-adjoint operator $-\A_\mu$. Define $a_\mu$ on the space $\D[a] = V_2^{1/2}(\R)$ by putting
\begin{equation*}
    a_\mu[u,v] = ((-\A)^{1/2} u, (-\A)^{1/2} v) - \mu \, u(a) \overline{v(a)} \text{,} \quad u, v \in \D[a] \text{.}
\end{equation*}
\begin{theorem}
    The form $a_\mu$ is semi-bounded from below and closed.
\end{theorem}
\begin{proof}
    Let $u \in \D[a]$. Using the Lemma \ref{l21} with $\kappa = \nu$, $\beta = 1/2$, we get
    \begin{equation*}
        a_\mu[u, u] = \| (-\A)^{1/2} u \|_2^2 - \mu \, |u(a)|^2 \geq - \nu \| u \|_2^2 \text{.}
    \end{equation*}
    
    The semi-boundedness from below is proved, let us prove closedness. We have
    \begin{equation*}
        \| u \|_a^2
        = a_\mu[u,u] + (\nu + 1) \| u \|_2^2
        = \| (-\A)^{1/2} u \|_2^2 - \mu |u(a)|^2 + (\nu + 1) \| u \|_2^2
    \end{equation*}
    \begin{equation*}
        \leq (\nu + 1) \| (-\A)^{1/2} u \|_2^2 + (\nu + 1) \| u \|_2^2
    \end{equation*}
    \begin{equation*}
        = \frac{\nu + 1}{2\pi} \int \limits_\R \big( 1 + \Psi(p) \big) |\widehat{u}(p)|^2 \, dp
        = \frac{\nu + 1}{2\pi} |u|_{1/2}^2 \text{,}
    \end{equation*}
    which means that any sequence converging in $V_{1/2}^2(\R)$ converges with respect to the norm $\|\cdot\|_a$. Thus, closedness has been proven.
\end{proof}

\begin{theorem}
    The operator $-\A_\mu$ corresponds to the form $a_\mu$.
\end{theorem}
\begin{proof}
    Let us show that the conditions \eqref{eq:op_and_form_cond1}, \eqref{eq:op_and_form_cond2} are met. 
    
    Let $u = \varphi + C \psi_\nu(\cdot - a) \in \D(\A_\mu) = \D(-\A_\mu)$. We have
    \begin{equation*}
        | u |_{1/2}^2
        = \int \limits_\R \big( 1 + \Psi(p) \big) |\widehat{u}(p)|^2 \, dp
    \end{equation*}
    \begin{equation*}
        \leq 2 \int \limits_\R \big( 1 + \Psi(p) \big) |\widehat{\varphi}(p)|^2 \, dp
        + 2 C^2 \int \limits_\R \big( 1 + \Psi(p) \big) |\widehat{\psi}_\nu(p)|^2 \, dp
    \end{equation*}
    \begin{equation*}
        \leq 2 \, | \varphi |_{1/2}^2
        + 2 C^2 \int \limits_\R \frac{\Psi(p) + 1}{(\Psi(p) + \nu)^2} \, dp < \infty \text{,}
    \end{equation*}
    therefore, $u \in \D(a_\mu)$, and the condition \eqref{eq:op_and_form_cond1} is satisfied.

    Let $u \in \D(-\A_\mu)$, $v \in \D(a_\mu)$. Using the Lemma \ref{l32}, we obtain
    \begin{equation*}
        (-\A_\mu u,v)
        = -\frac{1}{2\pi} (\widehat{\A_\mu u}, \widehat{v})
        = -\frac{1}{2\pi} (\widehat{\A}_\mu \widehat{u}, \widehat{v})
    \end{equation*}
    \begin{equation*}
        = \frac{1}{2\pi} \big( \Psi \widehat{u}, \widehat{v} \big)
        - \frac{1}{2\pi} \mu \, u(a) (e^{i(\cdot)a}, \widehat{v})
        = \frac{1}{2\pi} \big( \sqrt{\Psi} \widehat{u}, \sqrt{\Psi} \widehat{v} \big)
        - \mu \, u(a) \overline{v(a)}
    \end{equation*}
    \begin{equation*}
        = \big( (-\A)^{1/2} u, (-\A)^{1/2} v \big)
        - \mu \, u(a) \overline{v(a)} \text{,}
    \end{equation*}
    where, as earlier, $(e^{i(\cdot)x_k}, \widehat{v})$ is to be considered as
    \begin{equation*}
         \lim_{M \to \infty} \int \limits_{-M}^M e^{ipa} \, \overline{\widehat{v}(p)} \, dp \text{.}
    \end{equation*}
    
    Thus, the condition \eqref{eq:op_and_form_cond2} is also met, and the proof is complete.
\end{proof}

Let's describe the spectrum of the operator $\A_\mu$.

\begin{lemma}
    Let $\lambda \in \mathbb{C} \setminus \big((-\infty, 0] \cup \{ \nu \} \big)$. The resolvent of the operator $\A_\mu$ acts on $f\in L_2(\R)$ by the formula
    \begin{equation}
    \label{eq:resolvent}
        (\A_\mu - \lambda)^{-1} f
        = (\A - \lambda)^{-1} f + \frac{1}{\nu - \lambda} \, \frac{(f, \psi_{\bar{\lambda}}(\cdot - a))}{(\psi_\nu, \psi_{\bar{\lambda}})} \, \psi_{\lambda}(\cdot - a) \text{.}
    \end{equation}
\end{lemma}
\begin{proof}
    To obtain the formula, one can use the considerations on operators with one-rank perturbations given in \cite[ch. 11, 11.2]{Simon}.
    
    First, we show that the operator in the right part of \eqref{eq:resolvent} is bounded. Let $f \in L_2(\R)$. We have
    \begin{equation*}
        \| (\A_\mu - \lambda)^{-1} f \|_2^2
        \leq \| (\A - \lambda)^{-1} f \|_2^2
        + \frac{|(f, \psi_{\bar{\lambda}}(\cdot - a))|}{|(\nu - \lambda)(\psi_\nu, \psi_{\bar{\lambda}})|} \| \psi_\lambda(\cdot - a) \|_2^2
    \end{equation*}
    \begin{equation*}
        = \frac{1}{2 \pi} \bigg\| \frac{\widehat{f}}{\Psi + \lambda} \bigg\|_2^2
        + \frac{\| f \|_2^2 \| \psi_\lambda \|_2^2}{|(\nu - \lambda)(\psi_\nu, \psi_{\bar{\lambda}})|} \| \psi_\lambda \|_2^2
        \leq \Bigg( \frac{1}{|\lambda|} + \frac{\| \psi_\lambda \|_2^4}{|(\nu - \lambda)(\psi_\nu, \psi_{\bar{\lambda}})|} \Bigg) \| f \|_2^2 \text{.}
    \end{equation*}
    
    Now it is enough to check the formula for the functions from the range of $\A_\mu - \lambda$. Let
    \begin{equation*}
        f = (\A_\mu - \lambda) (\varphi + C \psi_\nu(\cdot - a))
        = (\A - \lambda) \varphi + C (\nu - \lambda) \psi_\nu(\cdot - a)
    \end{equation*}
    for some $\varphi \in V_2^1(\R) : \varphi(a) = 0$, $C \in \R$. Then
    \begin{equation*}
        \big( \F \big( (\A - \lambda)^{-1} f + \frac{1}{\nu - \lambda} \, \frac{(f, \psi_{\bar{\lambda}}(\cdot - a))}{(\psi_\nu, \psi_{\bar{\lambda}})} \, \psi_\lambda(\cdot - a) \big) \big)(p)
    \end{equation*}
    \begin{equation*}
        = - \frac{\widehat{f}(p)}{\Psi(p) + \lambda}
        + \frac{1}{2 \pi} \frac{1}{(\nu - \lambda)(\psi_\nu, \psi_{\bar{\lambda}})} \bigg( \widehat{f}, \frac{e^{i(\cdot)a}}{\Psi + \bar{\lambda}} \bigg) \frac{e^{ipa}}{\Psi(p) + \lambda}
    \end{equation*}
    \begin{equation*}
        = \widehat{\varphi}(p)
        - \frac{C (\nu - \lambda) e^{ipa}}{(\Psi(p) + \nu)(\Psi(p) + \lambda)}
    \end{equation*}
    \begin{equation*}
        + \frac{1}{(\nu - \lambda)(\psi_\nu, \psi_{\bar{\lambda}})} \bigg( -\frac{1}{2 \pi} \int \limits_\R e^{-ipa} \widehat{\varphi}(p) \, dp \bigg) \frac{e^{ipa}}{\Psi(p) + \lambda}
    \end{equation*}
    \begin{equation*}
        + \frac{1}{(\nu - \lambda)(\psi_\nu, \psi_{\bar{\lambda}})} \bigg( \frac{1}{2 \pi} \int \limits_\R \frac{C(\nu - \lambda)}{(\Psi(p) + \nu)(\Psi(p) + \lambda)} \, dp \bigg) \frac{e^{ipa}}{\Psi(p) + \lambda}
    \end{equation*}
    \begin{equation*}
        = \widehat{\varphi}(p)
        + \frac{C e^{ipa}}{\Psi(p) + \nu}
        - \frac{C e^{ipa}}{\Psi(p) + \lambda}
        - \frac{\varphi(a)}{(\nu - \lambda)(\psi_\nu, \psi_{\bar{\lambda}})} 
        + \frac{C e^{ipa}}{\Psi(p) + \lambda}
    \end{equation*}
    \begin{equation*}
        = \widehat{\varphi}(p)
        + C \widehat{\psi}_\nu(p) e^{ipa} \text{,}
    \end{equation*}
    which completes the proof.
\end{proof}
\begin{corollary*}
    The kernel $r_\mu(\lambda,x,y)$ of $(\A_\mu - \lambda)^{-1}$ is given by
    \begin{equation}
    \label{eq:res_kernel}
        r_\mu(\lambda,x,y)
        = - \frac{1}{2 \pi} \int \limits_\R \frac{e^{-ip(x-y)}}{\Psi(p) + \lambda} \, dp
        + \frac{\psi_\lambda(x - a) \psi_{\lambda}(y - a)}{(\nu - \lambda)(\psi_\nu, \psi_{\bar{\lambda}})} \text{.}
    \end{equation}
\end{corollary*}

\begin{theorem}
    The spectrum of the operator $\A_\mu$ consists of $(-\infty, 0]$, which is the continuous part, and a single eigenvalue -- $\nu$, which is the discrete part.
\end{theorem}
\begin{proof}
    The fact that $\nu$ is the eigenvalue of the operator $\A_\mu$ is evident both from the definition of $\A_\mu$ and from the formula for the resolvent \eqref{eq:resolvent}. From the same formula it follows that $\A_\mu$ inherits the spectrum of the operator $\A$, which is continuous and lies on $(-\infty, 0]$.
\end{proof}

%% file: 5-semigroups.tex
\section{Semigroups of operators $\{ U_t \}$, $\{ \widetilde{U}_t \}$}

Denote by $\{\F_t\}_{t\geq 0}$ the filtration generated by the process $\xi(t)$.

As shown in the previous section, $-\A_\mu$ is a semi-bounded from below self-adjoint operator. Spectral theory allows us to define \cite[ch. 8, \S 8.2]{BirSol} a semigroup of operators $\{ e^{t \A_\mu} \}_{t \geq 0}$ such that, given $f \in L_2(\R)$, the function
\begin{equation*}
    u(t,x) = (e^{t \A_\mu} f)(x)
\end{equation*}
is the unique solution to the Cauchy problem
\begin{equation}
\label{eq:cauchy_problem2}
    \frac{\partial u}{\partial t}
    = \A_\mu u \text{,} \quad u(0,x) = f(x) \text{,}
\end{equation}
in the domain
\begin{align*}
    \{ u: \R_+ \times \R \rightarrow \R \, \big| \, t \mapsto u(t, \cdot) \in C(\R_+, L_2(\R)) \text{,} \ &u(\cdot, x)\in C^1(\R) \text{,} \\ &u(t,\cdot) \in \D(\A_\mu) \} \text{.}
\end{align*}

Let us show that \eqref{eq:feynman_kac_formula} is the probabilistic representation of $\{ e^{t\A_\mu} \}$. First, we need the following statement.

\begin{lemma}
\label{l41}
    \item[1.] Let $f \in C_b(\R)$. Then
    \begin{equation*}
        \E f(\xi_x(t)) L(t,x - a)
        = \frac{1}{2 \pi} \int \limits_\R e^{-ip(x - a)} \bigg( \int \limits_0^t  e^{-\tau \Psi(p)} \, \E f(a - \xi(t - \tau)) \, d\tau \bigg) \, dp \text{.}
    \end{equation*}
    \item[2.] Let $f \in C_b(\R)$, $k \in \mathbb{N}$. Then
    \begin{equation*}
        \big| \E f(\xi_x(t)) (L(t,x -a))^k \big|
        \leq \frac{k!}{(2 \pi)^k} \, \E |f(\xi_x(t))| \, \Big( \int \limits_\R \frac{1 - e^{-t \Psi(p)}}{\Psi(p)} \, dp \Big)^k \text{.}
    \end{equation*}
\end{lemma}
\begin{proof}
    Let $f \in C_b(\R)$.
    We begin by proving the first part of the statement. We ger
    \begin{equation*}
        \E f(\xi_x(t)) L(t,x - a)
    \end{equation*}
    \begin{equation*}
        = \E f(\xi_x(t)) \Big( \lim_{M \to \infty} \frac{1}{2 \pi} \int \limits_{-M}^M e^{-iq(x - a)} \Big( \int \limits_0^t e^{iq\xi(\tau)} \, d\tau \Big) \, dq \Big)
    \end{equation*}
    \begin{equation*}
        = \E \lim_{M \to \infty} \frac{1}{2 \pi} \int \limits_{-M}^M e^{-iq(x - a)} \Big( \int \limits_0^t e^{iq(\xi(\tau))} f(\xi_x(\tau) - \xi(t - \tau)) \, d\tau \Big) \, dq
    \end{equation*}
    \begin{equation*}
        = \int \limits_0^t \E \lim_{M \to \infty} \frac{1}{2 \pi} \int \limits_{-M}^M e^{-iq(x - a)} \Big( e^{iq\xi(\tau)} f(\xi_x(\tau) - \xi(t - \tau)) \Big) \, dq \, d\tau
    \end{equation*}
    \begin{equation*}
        = \int \limits_0^t \E \Big( \E \lim_{M \to \infty} \frac{1}{2 \pi} \int \limits_{-M}^M e^{-iq(x - a)} \Big( e^{iq\xi(\tau)} f(\xi_x(\tau) - y) \Big) \, dq \, \big|_{y = \xi(t - \tau)} \Big) \, d\tau
    \end{equation*}
    \begin{equation*}
        = \int \limits_0^t \E \Big( \lim_{M \to \infty} \frac{1}{2 \pi} \int \limits_{-M}^M e^{-iq(x - a)} \Big( \E e^{iq\xi(\tau)} f(\xi_x(\tau) - y) \Big) \, dq \, \big|_{y = \xi(t - \tau)} \Big) \, d\tau
    \end{equation*}
    \begin{equation*}
        = \int \limits_0^t \E \lim_{M \to \infty} \int \limits_\R \delta_M(x - a - z) f(x - z - y) \Big( \int \limits_\R e^{-ipz} e^{-\tau \Psi(p)} \, dp \Big) \, dz \, \big|_{y = \xi(t - \tau)} \, d\tau \text{,}
    \end{equation*}
    where
    \begin{equation*}
        \delta_M(u)
        = \Big( \frac{1}{2 \pi} \Big)^2 \int \limits_{-M}^M e^{-iqu} \, dq
        = \frac{1}{2\pi} \frac{\text{sin}(Mu)}{\pi u} \text{.}
    \end{equation*}
    
    Furthermore,
    \begin{equation*}
        \int \limits_0^t \E \lim_{M \to \infty} \int \limits_\R q_M(x - a - z) f(x - z - y) \Big( \int \limits_\R e^{-ipz} e^{-\tau \Psi(p)} \, dp \Big) \, dz \, \big|_{y = \xi(t - \tau)} \, d\tau \text{,}
    \end{equation*}
    \begin{equation*}
        = \frac{1}{2 \pi} \int \limits_0^t \E f(a - \xi(t - \tau)) \Big( \int \limits_\R e^{-ip(x - a)} e^{-\tau \Psi(p)} \, dp \Big) \, d\tau
    \end{equation*}
    \begin{equation*}
        = \frac{1}{2 \pi} \int \limits_\R e^{-ip(x - a)} \Big( \int \limits_0^t e^{-\tau \Psi(p)} \, \E f(a - \xi(t - \tau)) \, d\tau \Big) \, dp \text{.}
    \end{equation*}

    We proceed with the second part of the statement. Let $k \in \mathbb{N}$ and $\Theta_t^k = \{ (\tau_1, \ldots, \tau_k) \, |\, 0 \leq\tau_1 \leq\ldots \leq\tau_k \leq t \}$. Then
    \begin{equation*}
        \big| \E f(\xi_x(t)) (L(t,x - a))^k \big|
    \end{equation*}
    \begin{equation*}
        = \bigg| \E f(\xi_x(t)) \bigg( \frac{1}{2 \pi}  \int \limits_\R e^{-ip(x - a)} \Big( \int \limits_0^t e^{ip\xi(\tau)} \, d\tau \Big) \, dp \bigg)^k \bigg|
    \end{equation*}
    \begin{equation*}
        \leq \Big( \frac{1}{2 \pi} \Big)^k \int \limits_{\R^k} \int \limits_{[0,t]^k} \big| \E f(\xi_x(t)) \, \prod_{l=1}^k e^{ip_l\xi(\tau_l)} \big| \, d\tau \, dp
    \end{equation*}
    \begin{equation*}
        \leq \frac{k!}{(2 \pi)^k} \, \int \limits_{\R^k} \int \limits_{\Theta_t^k} \big| \E f(\xi_x(t)) \, \prod_{l=1}^k e^{ip_l\xi(\tau_l)} \big| \, d\tau \, dp \text{.}
    \end{equation*}
    
    Here and subsequently, we consider the product over the empty set of indices to be zero. Let's make a change introducing the variables $q_1 \text{,}\, \ldots\text{, } q_k$, where
    \begin{equation*}
        q_l = \sum_{m=l}^k p_m \text{,} \quad 1 \leq l \leq k \text{,}
    \end{equation*}
    and put $\tau_0 = 0$. We obtain
    \begin{equation*}
        \frac{k!}{(2 \pi)^k} \, \int \limits_{\R^k} \int \limits_{\Theta_t^k} \big| \E f(\xi_x(t)) \, \prod_{l=1}^k e^{ip_l\xi(\tau_l)} \big| \, d\tau \, dp
    \end{equation*}
    \begin{equation*}
        = \frac{k!}{(2 \pi)^k} \, \int \limits_{\R^k} \int \limits_{\Theta_t^k} \big| \E f(\xi_x(t)) \prod_{l=1}^k e^{i q_l (\xi(\tau_l) - \xi(\tau_{l-1}))} \big| \, d\tau \, dq
    \end{equation*}
    \begin{equation*}
        = \frac{k!}{(2 \pi)^k} \,\int \limits_{\R^k} \int \limits_{\Theta_t^k} \Big| \E \Big( \E f(y - \xi(t - \tau_l)) \prod_{l=1}^k e^{i q_l (\xi(\tau_l) - \xi(\tau_{l-1}))} \, \big|_{y = \xi_x(\tau_l)} \Big) \Big| \, d\tau \, dq
    \end{equation*}
    \begin{equation*}
        = \frac{k!}{(2 \pi)^k} \,\int \limits_{\R^k} \int \limits_{\Theta_t^k} \Big| \E \Big( \E f(y - \xi(t - \tau_l)) \prod_{l=1}^k \E e^{i q_l (\xi(\tau_l) - \xi(\tau_{l-1}))} \, \big|_{y = \xi_x(\tau_l)} \Big) \Big| \, d\tau \, dq
    \end{equation*}
    \begin{equation*}
        = \frac{k!}{(2 \pi)^k} \,\int \limits_{\R^k} \int \limits_{\Theta_t^k} \Big| \E \Big( \E f(y - \xi(t - \tau_l)) \, \big|_{y = \xi_x(\tau_l)} \Big) \prod_{l=1}^k e^{-(\tau_l - \tau_{l-1}) \Psi(q_l)} \Big| \, d\tau \, dq
    \end{equation*}
    \begin{equation*}
        \leq \frac{k!}{(2 \pi)^k} \int \limits_{\R^k} \int \limits_{\Theta_t^k} \E |f(\xi_x(t))| \prod_{l=1}^k e^{-(\tau_l - \tau_{l-1}) \Psi(q_l)} \, d\tau \, dq \text{.}
    \end{equation*}
    
    Define $\Xi_t^k = \{ (s_1, \ldots, s_k) \, | \, s_1 + \ldots + s_k \leq t\}$ and make a change introducing the variables $s_1 \text{, } \ldots \text{,} s_k$, where
    \begin{equation*}
        s_1 = \tau_1 \text{,} \quad s_l = \tau_l - \tau_{l-1} \text{,} \quad 2 \leq l \leq k \text{.}
    \end{equation*}
    We obtain
    \begin{equation*}
        \frac{k!}{(2 \pi)^k} \int \limits_{\R^k} \int \limits_{\Omega_t^k} \E |f(\xi_x(t))| \, e^{-\tau_1 \Psi(q_1)} \prod_{l=2}^k e^{-(\tau_l - \tau_{l-1}) \Psi(q_l)} \, d\tau \, dq
    \end{equation*}
    \begin{equation*}
        = \frac{k!}{(2 \pi)^k} \, \E |f(\xi_x(t))| \int \limits_{\R^k} \int \limits_{\Xi_t^k} \, \prod_{l=1}^k e^{-s_l \Psi(q_l)} \, ds \, dq
    \end{equation*}
    \begin{equation*}
        \leq \frac{k!}{(2 \pi)^k} \, \E |f(\xi_x(t))| \int \limits_{\R^k} \int \limits_{[0,t]^k} \, \prod_{l=1}^k e^{-s_l \Psi(q_l)} \, ds \, dq
    \end{equation*}
    \begin{equation*}
        = \frac{k!}{(2 \pi)^k} \, \E |f(\xi_x(t))| \, \int \limits_{\R^k} \, \prod_{l=1}^k \frac{1 - e^{-t \Psi(q_l)}}{\Psi(q_l)} \, dq
    \end{equation*}
    \begin{equation*}
        = \frac{k!}{(2 \pi)^k} \, \E |f(\xi_x(t))| \, \bigg( \int \limits_\R \frac{1 - e^{-t \Psi(p)}}{\Psi(p)} \, dp \bigg)^k \text{.}
    \end{equation*}
\end{proof}

Let $\{ U_t \}_{t \geq 0}$ be a family of operators acting on $C_b(\R)$ by the formula
\begin{equation*}
    (U_t f)(x)
    = \E f(\xi_x(t)) e^{\mu L(t,x - a)} \text{,} \quad f \in C_b(\R) \text{.}
\end{equation*}

\begin{theorem}
\label{th51}
    \item[1.] The family $\{ U_t \}$ is a $C_0$-semigroup in $L_2(\R)$ with the generator $\A_\mu$.
    \item[2.] For any $f \in L_2(\R) \cap C_b(\R)$
    \begin{equation*}
        (e^{t \A_\mu} f)(x)
        = \E f(\xi_x(t)) e^{\mu L(t,x - a)} \text{.}
    \end{equation*}
\end{theorem}
\begin{proof}
    Let $f \in L_2(\R) \cap C_b(\R)$. 
    
    Obviously,
    \begin{equation*}
        (U_0 f)(x) = f(x) \text{.}
    \end{equation*}
    
    Let's check that the semigroup property is met. We have
    \begin{equation*}
        (U_s (U_t f))(x)
        = \E \Big( e^{\mu L(s,x - a)} \Big[ \E \big(  f(\xi_y(t)) e^{\mu L(t,y - a)} \big) \big|_{\xi_x(x) = y} \Big] \Big)
    \end{equation*}
    \begin{equation*}
        = \E \big( e^{\mu L(s,x - a)} \, \E \big( f(\xi_{\xi_x(s)}(t)) e^{\mu L(t,\xi_x(s) - a)} \big| \F_s \big) \big)
    \end{equation*}
    \begin{equation*}
        = \E \, \E \big( f(\xi_{\xi_x(s)}(t)) e^{\mu L(s,x - a)} e^{\mu L(t,\xi_x(s) - a)} \big| \F_s \big)
    \end{equation*}
    \begin{equation*}
        = \E f(\xi_{\xi_x(s)}(t)) e^{\mu L(s,x - a)} e^{\mu L(t,\xi_x(s) - a)}
    \end{equation*}
    \begin{equation*}
        = \E f(\xi_x(t+s)) e^{\mu L(t+s,x - a)}
        = (U_{t+s} f)(x) \text{.}
    \end{equation*}
    
    Let's proceed with strong continuity. We have
    \begin{equation*}
        \| U_t f  - f \|_2
        = \| (T_t f - f) + (U_t f - T_t f) \|_2
    \end{equation*}
    \begin{equation}
    \label{th51_proof1}
        \leq \| T_t f - f \|_2
        + \| U_t f - T_t f \|_2 \text{,}
    \end{equation}
    where $\{ T_t \}_{t \geq 0}$ is the semigroup generated by $\xi(t)$.

    The first term in \eqref{th51_proof1} tends to zero as $t \to 0+$, since $f$ belongs to the domain of the $C_0$-semigroup $\{T_t\}$.

    Consider the second term in \eqref{th51_proof1}. Using the second statement of Lemma \ref{l41}, we obtain
    \begin{equation*}
        \| U_t f - T_t f \|_2
        = \| \E f(\xi_x(t)) (e^{\mu L(t,x - a)} - 1) \|_2
    \end{equation*}
    \begin{equation*}
        \leq \sum_{k=1}^\infty \frac{\mu^k}{k!} \| \E f(\xi_x(t)) (L(t,x))^k \|_2
    \end{equation*}
    \begin{equation*}
        \leq \sum_{k=1}^\infty \frac{\mu^k}{(2 \pi)^k} \bigg \| \E |f(\xi_x(t))| \, \Big( \int \limits_\R \frac{1 - e^{-t \Psi(p)}}{\Psi(p)} \, dp \Big)^k \bigg \|_2
    \end{equation*}
    \begin{equation*}
        = \sum_{k=1}^\infty \frac{\mu^k}{(2 \pi)^k} \Big( \int \limits_\R \frac{1 - e^{-t \Psi(p)}}{\Psi(p)} \, dp \Big)^k \| f \|_2 \stackrel[t \to 0+]{}{\longrightarrow} 0 \text{.}
    \end{equation*}

    Due to the density of $L_2(\R) \cap C_b(\R)$ in $L_2(\R)$ the semigroup properties are also met for the functions from $L_2(\R)$.
    
    Now, let us evaluate the generator. Let $f \in \D(\A_\mu)$, $f = \varphi + C \psi_\nu(\cdot - a)$. Then
    \begin{equation*}
        \Big \| \frac{U_t f - f}{t} - \A_\mu f \Big \|_2
    \end{equation*}
    \begin{equation}
    \label{th51_proof2}
        \leq \Big \| \frac{U_t \varphi - \varphi}{t} - \A \varphi \Big \|_2
        + \Big \| \frac{U_t \psi_\nu(\cdot - a) - \psi_\nu(\cdot - a)}{t} - \nu \psi_\nu(\cdot - a) \Big \|_2 \text{.}
    \end{equation}

    By $I_1$ and $I_2$ respectively, denote the terms in \eqref{th51_proof2}. For $I_1$, we have
    \begin{equation*}
        I_1
        = \Big \| \frac{U_t \varphi - \varphi}{t} - \A \varphi \Big \|_2
        = \Big \| \frac{T_t \varphi - \varphi}{t} - \A \varphi + \frac{U_t \varphi - T_t \varphi}{t} \Big \|_2
    \end{equation*}
    \begin{equation*}
        \leq \Big \| \frac{T_t \varphi - \varphi}{t} - \A \varphi \Big \|_2 + \Big \| \frac{U_t \varphi - T_t \varphi}{t} \Big \|_2 \text{.}
    \end{equation*}
    The first term tends to zero as $t \to 0+$, since $\A$ is the generator of $\{T_t\}$. Evaluating the second term using the second statement of Lemma \ref{l41}, we obtain
    \begin{equation*}
        \Big \| \frac{U_t \varphi - T_t \varphi}{t} \Big \|_2
        = \Big \| \frac{1}{t} \, \E \varphi(\cdot - \xi(t)) (e^{\mu L(t,\cdot - a)} - 1) \Big \|_2
    \end{equation*}
    \begin{equation*}
        \leq \Big \| \frac{\mu}{t} \, \E \varphi(\cdot - \xi(t)) L(t,\cdot - a) \Big \|_2
        + \sum_{k=2}^\infty \frac{\mu^k}{(2 \pi)^k \, t} \Big( \int \limits_\R \frac{1 - e^{-t \Psi(p)}}{\Psi(p)} \, dp \Big)^k \| \varphi \|_2 \text{,}
    \end{equation*}
    where the second term tends to zero as $t \to 0+$.
    
    Furthermore, from the first statement of Lemma \ref{l41} it follows that
    \begin{equation*}
        \bigg \| \frac{\mu}{t} \, \E \varphi(\cdot - \xi(t)) L(t,\cdot - a) \bigg \|_2
    \end{equation*}
    \begin{equation*}
        = \bigg \| \frac{1}{2 \pi} \int \limits_\R e^{-ip(\cdot - a)} \bigg( \frac{\mu}{t} \int \limits_0^t e^{-\tau \Psi(p)} \, \E \varphi(a - \xi(t - \tau)) \, d\tau \bigg) \, dp \bigg \|_2
    \end{equation*}
    \begin{equation*}
        = \frac{1}{\sqrt{2 \pi}} \bigg \| \frac{\mu}{t} \, \int \limits_0^t e^{-\tau \Psi(\cdot)} \, \E \varphi(a - \xi(t - \tau)) \, d\tau \bigg \|_2 \text{.}
    \end{equation*}
    
    Using the fact that $T_s \varphi = e^{s \A} \varphi$ and $\varphi(a) = 0$, we get
    \begin{equation*}
        \frac{1}{\sqrt{2 \pi}} \bigg \| \frac{\mu}{t} \, \int \limits_0^t e^{-\tau \Psi(\cdot)} \, \E \varphi(a - \xi(t - \tau)) \, d\tau \bigg \|_2
    \end{equation*}
    \begin{equation*}
        = \frac{1}{\sqrt{2 \pi}} \bigg \| \frac{\mu}{t} \, \int \limits_0^t e^{-\tau \Psi(\cdot)} \, \bigg( \frac{1}{2 \pi} \int \limits_\R e^{-iqa} \big( e^{-(t - \tau) \Psi(q)} - 1 \big) \widehat{\varphi}(q) \, dq \bigg) \, d\tau \bigg \|_2
    \end{equation*}
    \begin{equation*}
        \leq \frac{1}{\sqrt{2 \pi}} \bigg \| \frac{\mu}{t} \, \int \limits_0^t e^{-\tau \Psi(\cdot)} \, \bigg( \frac{1}{2 \pi} \int \limits_\R \big( 1 - e^{-t \Psi(q)} \big) |\widehat{\varphi}(q)| \, dq \bigg) \, d\tau \bigg \|_2
    \end{equation*}
    \begin{equation*}
        = \frac{1}{\sqrt{2 \pi}} \bigg \| \frac{\mu}{t} \int \limits_0^t e^{-\tau \Psi(\cdot)} \Big( \frac{1}{2 \pi} \int \limits_\R \big( 1 - e^{-t \Psi(q)} \big) \frac{\sqrt{1 + \Psi^2(q)}}{\sqrt{1 + \Psi^2(q)}} \, |\widehat{\varphi}(q)| \, dq \Big) \, d\tau \bigg \|_2
    \end{equation*}
    \begin{equation*}
        = \frac{\mu}{(2 \pi)^{3/2}} \bigg( \int \limits_\R \frac{(1 - e^{-t \Psi(q)})^2}{t \, (1 + \Psi^2(q))} \, dq \bigg)^{1/2} \, \bigg( \int \limits_\R \bigg( \frac{1 - e^{-t \Psi(p)}}{\Psi(p)} \bigg)^2 \, dp \bigg)^{1/2} |\varphi|_1 \stackrel[t \to 0+]{}{\longrightarrow} 0 \text{.}
    \end{equation*}

    Now, let's evaluate $I_2$. To do this, we need a statement from \cite{SalYor}. Introduce it in our notation.
    \begin{theorem}[{[Salminen, Yor, 2007]}]
        There exists a $\F_t$-martingale with zero mean $M_{\nu, x - a}(t)$ such that
        \begin{equation}
        \label{th51_proof3}
            \psi_\nu(\xi_x(t)- a) - \psi_\nu(x - a)
            = \nu \int_0^t \psi_\nu(\xi_x(\tau) - a) \, d\tau - L(t,x - a) + M_{\nu,x - a}(t) \text{.}
        \end{equation}
    \end{theorem}

    Using the formula \eqref{th51_proof3} and the second statement of Lemma \ref{l41}, we obtain
    \begin{equation*}
        I_2
        = \Big \| \frac{U_t \psi_\nu(\cdot - a) - \psi_\nu(\cdot - a)}{t} - \nu \psi_\nu(\cdot - a) \Big \|_2
    \end{equation*}
    \begin{equation*}
        = \bigg \| \frac{1}{t} U_t \psi_\nu(\cdot - a) + \frac{1}{t} \Big( \E \psi_\nu(\cdot - \xi(t) - a) - \psi_\nu(\cdot - a) \Big) - \nu \psi_\nu(\cdot - a) \bigg \|_2
    \end{equation*}
    \begin{equation*}
        = \bigg \| \frac{1}{t} U_t \psi_\nu(\cdot - a) + \frac{1}{t} \Big( \nu \E \int \limits_0^t \psi_\nu(\cdot - \xi(\tau) - a) \, d\tau - L(t,\cdot - a) \Big) - \nu \psi_\nu(\cdot - a) \bigg \|_2
    \end{equation*}
    \begin{equation}
    \label{th51_proof4}
        \leq \Big \| \frac{1}{t} \Big( \mu \, \E \psi_\nu(\cdot - \xi(t) - a) - 1 \Big) L(t,\cdot - a) \Big \|_2
    \end{equation}
    \begin{equation}
    \label{th51_proof5}
        + \, \Big \| \nu \Big( \frac{1}{t} \, \E \int \limits_0^t \psi_\nu(\cdot - \xi(\tau) - a) \, d\tau - \psi_\nu(\cdot - a) \Big) \Big \|_2
    \end{equation}
    \begin{equation*}    
        + \, \sum_{k=2}^\infty \frac{\mu^k}{(2 \pi)^k \, t} \Big( \int \limits_\R \frac{1 - e^{-t \Psi(p)}}{\Psi(p)} \, dp \Big)^k \| \psi_\nu \|_2 \text{.}
    \end{equation*}

    The last term tends to zero as $t \to 0+$. Let us show that the same is fair for the first two terms. By $J_1$ and $J_2$ respectively, denote \eqref{th51_proof4} and \eqref{th51_proof5}. Using the first statement of Lemma \ref{l41} and the fact that $\mu\, \psi_\nu(a) =1$, we obtain
    \begin{equation*}
        J_1 = \frac{1}{\sqrt{2 \pi}} \bigg \| \frac{1}{t} \int \limits_0^t e^{-\tau \Psi(\cdot)} \big( \mu \E \psi_\nu(a - \xi(t - \tau)) - 1 \big) \, d\tau \bigg \|_2
    \end{equation*}
    \begin{equation*}
        = \frac{1}{\sqrt{2 \pi}} \bigg \| \frac{\mu}{t} \int \limits_0^t e^{-\tau \Psi(\cdot)} \big( \E \psi_\nu(a - \xi(t - \tau)) - \psi_\nu(a) \big) \, d\tau \bigg \|_2
    \end{equation*}
    \begin{equation*}
        = \frac{1}{\sqrt{2 \pi}} \bigg \| \frac{\mu}{t} \int \limits_0^t e^{-\tau \Psi(\cdot)} \Big( \frac{1}{2 \pi} \int \limits_\R e^{-iqa} \big( e^{-(t - \tau) \Psi(q)} - 1 \big) \widehat{\psi}_\nu(q) \, dq \Big) \, d\tau \bigg \|_2
    \end{equation*}
    \begin{equation*}
        \leq \frac{1}{\sqrt{2 \pi}} \bigg \| \frac{\mu}{t} \int \limits_0^t e^{-\tau \Psi(\cdot)} \Big( \frac{1}{2 \pi} \int \limits_\R \big( 1 - e^{-t \Psi(q)} \big) \frac{1}{\sqrt{\Psi(q) + \nu}} \frac{1}{\sqrt{\Psi(q) + \nu}} \, dq \Big) \, d\tau \bigg \|_2
    \end{equation*}
    \begin{equation*}
        = \frac{1}{\sqrt{2 \pi}} \bigg( \int \limits_\R \frac{(1 - e^{-t \Psi(q)})^2}{t \, (\Psi(q) + \nu)} \, dq \bigg)^{1/2} \bigg( \int \limits_\R \bigg( \frac{1 - e^{-t \Psi(p)}}{\Psi(p)} \bigg)^2 \, dp \bigg)^{1/2} \stackrel[t \to 0+]{}{\longrightarrow} 0 \text{.}
    \end{equation*}

    Now, let's evaluate $J_2$. Using definition of the local time, we have
    \begin{equation*}
        \Big \| \nu \Big( \frac{1}{t} \, \E \int_0^t \psi_\nu(\cdot - \xi(\tau) - a) \, d\tau - \psi_\nu(\cdot - a) \Big) \Big \|_2
    \end{equation*}
    \begin{equation*}
        = \Big \| \nu \Big( \frac{1}{t} \, \E \int \limits_\R \psi_\nu(\cdot - y - a) L(t, y) \, dy - \psi_\nu(\cdot - a) \Big) \Big \|_2
    \end{equation*}
    \begin{equation*}
        = \Big \| \nu \Big( \frac{1}{2\pi t} \, \E \int \limits_\R \frac{e^{-ip(\cdot - a)}}{\Psi(p) + \nu} \, \Big( \int \limits_0^t e^{ip\xi(\tau)} \, d\tau \Big) \, dp - \frac{1}{2 \pi} \int \limits_\R \frac{e^{-ip(\cdot - a)}}{\Psi(p) + \nu} \, dp \Big) \Big \|_2
    \end{equation*}
    \begin{equation*}
        = \Big \| \frac{\nu}{2\pi} \int \limits_\R \frac{e^{-ip(\cdot - a)}}{\Psi(p) + \nu} \, \Big( \frac{1}{t} \, \E \int \limits_0^t e^{ip\xi(\tau)} \, d\tau - 1 \Big) \, dp \Big \|_2
    \end{equation*}
    \begin{equation*}
        = \Big \| \frac{\nu}{2\pi} \int \limits_\R \frac{e^{-ip(\cdot - a)}}{\Psi(p) + \nu} \, \frac{1 - e^{-t \Psi(p)} - t \Psi(p)}{t \Psi(p)} \, dp \Big \|_2
    \end{equation*}
    \begin{equation*}
        = \frac{1}{\sqrt{2\pi}} \Big \| \frac{\nu}{\Psi + \nu} \, \frac{e^{-t \Psi} - 1 + t \Psi}{t \Psi} \Big \|_2 \stackrel[t \to 0+]{}{\longrightarrow} 0 \text{.}
    \end{equation*}
    
    Thus, $\A_\mu$ is the generator of $\{ U_t \}$ in $L_2(\R)$, which implies the second part of the statement of the theorem.
\end{proof}

Using the eigenfunction of $\A_\mu$, construct a $\F_t$-martingale. Define
\begin{equation*}
    \eta_\nu(t,x) = e^{-\nu t} \, \psi_\nu(\xi_x(t) - a) \, e^{\mu L(t,x - a)} \text{.}
\end{equation*}

\begin{theorem}
\label{th52}
    The process $\eta_\nu(t,x)$ is a $\F_t$-martingale.
\end{theorem}
\begin{proof}
    Let $\tau \leq t$. We have
    \begin{equation*}
        \E (\eta_\nu(t,x) \big| \F_\tau)
        = \E (e^{-\nu t} \, \psi_\nu(\xi_x(t) - a) \, e^{\mu L(t,x - a)} \big| \F_\tau)
    \end{equation*}
    \begin{equation*}
        = e^{-\nu t} \, \E (\psi_\nu(\xi_x(t) - a) \, e^{\mu \int \limits_0^t \, \delta(\xi_x(t) - a) d\tau} \big| \F_\tau)
    \end{equation*}
    \begin{equation*}
        = e^{-\nu t} \, e^{\mu \int \limits_0^\tau \, \delta(\xi_x(t) - a) d\tau} \, \E \big( \psi_\nu \big( \xi_x(\tau) - (\xi_x(t) - \xi_x(\tau)) - a \big) \, e^{\mu \int \limits_\tau^t \, \delta(\xi_x(t) - a) d\tau} \big| \F_\tau \big)
    \end{equation*}
    \begin{equation*}
        = e^{-\nu t} \, e^{\mu \int \limits_0^\tau \, \delta(\xi_x(t) - a) d\tau} \, \E \Big( \psi_\nu(\xi_y(t-\tau) - a) \, e^{\mu \int \limits_0^{t - \tau} \, \delta(\xi_y(t) - a) d\tau} \Big) \Big|_{y = \xi_x(\tau)} \text{.}
    \end{equation*}

    Since $\psi_\nu$ and $\nu$ respectively are the eigenfunction and the eigenvalue of the operator $\A_\mu$, $\psi_\nu$ and $e^{\nu \tau}$ respectively are the eigenfunction and the eigenvalue of the operator $e^{\tau \A_\mu}$. Thus, we have
    \begin{equation*}
        e^{-\nu t} \, e^{\mu \int \limits_0^\tau \, \delta(\xi_x(t) - a) d\tau} \, \E \Big( \psi_\nu(\xi_y(t-\tau) - a) \, e^{\mu \int \limits_0^{t - \tau} \, \delta(\xi_y(t) - a) d\tau} \Big) \Big|_{y = \xi_x(\tau)}
    \end{equation*}
    \begin{equation*}
        = e^{-\nu t} \, e^{\mu \int \limits_0^\tau \, \delta(\xi_x(t) - a) d\tau} e^{\nu (t - \tau)} \psi_\nu(\xi_x(\tau) - a)
    \end{equation*}
    \begin{equation*}
        = e^{-\nu \tau} \, \psi_\nu(\xi_x(\tau) - a) \, e^{\mu L(\tau,x - a)}
        = \eta_\nu(\tau,x) \text{,}
    \end{equation*}
    which completes the proof.
\end{proof}

The following statement is a consequence of martingality of the process $\eta_\nu(t,x)$.
\begin{theorem}
    For any $f \in L_2(\R)$
\begin{equation*}
        L_2 \, \text{-} \lim_{t \to \infty} e^{-\nu t} \, \E f(\xi_x(t)) e^{\mu L(t,x - a)} = \frac{(f, \psi_\nu(\cdot - a))}{\| \psi_\nu \|_2^2} \, \psi_\nu(x - a) \text{.}
    \end{equation*}
\end{theorem}
\begin{proof}
    Let $f \in L_2(\R)$. Use the decomposition
    \begin{equation*}
        f = f_0 + \frac{(f, \psi_\nu(\cdot - a))}{\| \psi_\nu \|_2^2} \, \psi_\nu(\cdot - a) \text{,}
    \end{equation*}
    where $f_0$ is orthogonal to $\psi_\nu(\cdot - a)$ in $L_2(\R)$.

    We have
    \begin{equation*}
        e^{-\nu t} \, \E f(\xi_x(t)) e^{\mu L(t,x - a)}
    \end{equation*}
    \begin{equation}
    \label{th53_proof1}
        = e^{-\nu t} \, \E f_0(\xi_x(t)) e^{\mu L(t,x - a)}
        + e^{-\nu t} \, \frac{(f, \psi_\nu(\cdot - a))}{\| \psi_\nu \|_2^2} \, \E \psi_\nu(\xi_x(t) - a) e^{\mu L(t,x - a)} \text{.}
    \end{equation}
    \begin{equation*}
        = e^{-\nu t} \, \E f_0(\xi_x(t)) e^{\mu L(t,x - a)}
        + \frac{(f, \psi_\nu(\cdot - a))}{\| \psi_\nu \|_2^2} \, \E \eta_\nu(t,x)
    \end{equation*}
    By Theorem \ref{th52}, the second term in \eqref{th53_proof1} is
    \begin{equation*}
        \frac{(f, \psi_\nu(\cdot - a))}{\| \psi_\nu \|_2^2} \, \psi_\nu(x - a) \text{.}
    \end{equation*}

    Show that the first term in \eqref{th53_proof1} tends to zero at $t \to \infty$. The function $f_0$ is orthogonal to $\psi_\nu(\cdot - a)$, which is the eigenfunction with the eigenvalue $\nu$, the only positive value of the spectrum of the operator $\A_\mu$. Therefore
    \begin{equation*}
        \| e^{-\nu t} \, \E f_0(\xi_x(t)) e^{\mu L(t,x - a)} \|_2
        = e^{-\nu t} \| e^{t \A_\mu} f_0 \|_2
    \end{equation*}
    \begin{equation*}
        \leq e^{-\nu t} \| f_0 \|_2
        \leq e^{-\nu t} \| f \|_2 \stackrel[t \to \infty]{}{\longrightarrow} 0 \text{.}
    \end{equation*}
\end{proof}

Now, we construct a Feller semigroup using $\eta_\nu(t,x)$.

Consider the space $C_0(\R)$ of the continuous, vanishing at infinity functions for which a functional $\|\cdot\|_\infty$ is finite, where
\begin{equation*}
    \| g \|_\infty = \sup_{x \in \R}{|g(x)|} \text{.}
\end{equation*}
The functional $\|\cdot\|_\infty$ is of course a norm, and $(C_0(\R), \|\cdot\|_\infty)$, or simply $C_0(\R)$, is a Banach space.

Let $\{ \widetilde{U}_t \}_{t \geq 0}$ be a family of operators acting on $C_0(\R)$ by the formula
\begin{equation*}
    ( \widetilde{U}_t g )(x)
    = \frac{\E \eta_\nu(t,x) g(\xi_x(t))}{\psi_\nu(x - a)} \text{,} \quad g \in C_0(\R) \text{.} 
\end{equation*}
\begin{theorem}
\label{th54}
    The family $\{\widetilde{U}_t\}$ is a Feller semigroup.
\begin{proof}
    Let $g\in C_0(\R)$. First, we prove that the family $\{\widetilde{U}_t\}$ is a semigroup.

    We have
    \begin{equation*}
        ( \widetilde{U}_0 g )(x)
        = \frac{\E \eta_\nu(0,x) g(\xi_x(0))}{\psi_\nu(x - a)}
        = g(x) \text{.}
    \end{equation*}

    Let $t \text{, } s\geq 0$. Then
    \begin{equation*}
        (\widetilde{U}_{t+s} f)(x)
        = \frac{e^{-(t + s)}}{\psi_\nu(x - a)} \big( U_{t+s} (\psi_\nu(\cdot - a) g) \big)(x)
    \end{equation*}
    \begin{equation*}
        = \frac{e^{-(t + s)}}{\psi_\nu(x - a)} \big( U_s ( U_t (\psi_\nu(\cdot - a) g)) \big)(x)
    \end{equation*}
    \begin{equation*}
        = \frac{e^{-s}}{\psi_\nu(x - a)} \Big( U_s \Big[ \psi_\nu(\cdot - a) \, \frac{e^{-t}}{\psi_\nu(\cdot - a)} \, U_t (\psi_\nu(\cdot - a) g) \Big] \Big)(x)
    \end{equation*}
    \begin{equation*}
        = \frac{e^{-s}}{\psi_\nu(x - a)} \big( U_s \big( \psi_\nu(\cdot - a) \, \widetilde{U}_t g \big) \big)(x)
        = (\widetilde{U}_s (\widetilde{U}_t f))(x) \text{.}
    \end{equation*}
    
    Now, let's prove Fellerness. Recall \cite[ch. III, 2.6]{RogWil} that for the semigroup $\{ \widetilde{U}_t \}$ to be a Feller one it suffices that for any $g \in C_0(\R)$
    \begin{itemize}
        \Item[a)]
        \begin{equation*}
            0 \leq g \leq 1
            \Rightarrow 0 \leq \widetilde{U}_t g \leq 1 \text{;}
        \end{equation*}
        \Item[b)]
        \begin{equation*}
            \lim_{t \to 0+} (\widetilde{U}_t g)(x)
            = g(x) \text{,} \quad x \in \R \text{.}
        \end{equation*}
    \end{itemize}
    
    If $0 \leq g \leq 1$, then by the Theorem \ref{th52}
    \begin{equation*}
        (\widetilde{U}_t g)(x)
        = \frac{\E \eta_\nu(t, x) g(\xi_x(t))}{\psi_\nu(x - a)}
        \leq \frac{\E \eta_\nu(t,x)}{\psi_\nu(x - a)}
        = 1 \text{.}
    \end{equation*}

    Furthermore, from continuity of $g$ and $\psi_\nu(\cdot - a)$ and the fact that the sample paths of the process $\xi(t)$ are right-continuous with probability 1, it follows that for any $x \in \R$ almost surely
    \begin{equation*}
        \frac{\eta_\nu(t,x) g(\xi_x(t))}{\psi_\nu(x - a)} \stackrel[t \to 0+]{}{\longrightarrow} g(x) \text{.}
    \end{equation*}
    Therefore
    \begin{equation*}
        \frac{\E \eta_\nu(t,x) g(\xi_x(t))}{\psi_\nu(x - a)} \stackrel[t \to 0+]{}{\longrightarrow} g(x) \text{.}
    \end{equation*}
    
    Thus, both conditions are met, hence $\{\widetilde{U}_t\}$ is a Feller semigroup.
\end{proof}
\end{theorem}

%% file: 6-attracting_distribution.tex
\section{Penalization}

In this section, we construct the measure \eqref{eq:penal_meas}. First, let's briefly describe the considerations we use.

Consider the measure \eqref{eq:gibbs_meas} on the sample paths of the process $\xi_x(t)$, $t \leq T$, in the case of a ``classical" potential $V$, and denote it by $\mathbf{Q}_{T,x}^V$. By $p_V(t,x,y)$ we denote the kernel of the operator $e^{t(\A + V)}$. The finite-dimensional distributions of $\mathbf{Q}_{T,x}^V$ are represented as follows.
\begin{equation*}
    \mathbf{Q}_{T,x}^V \{ \omega(t_1) \in B_1, \dots, \omega(t_n) \in B_n \}
\end{equation*}
\begin{equation*}
    = \frac{1}{Z_V(T,x)}
    \int \limits_{B_1} \ldots \int \limits_{B_n} p_V(t_1, x, x_1) \prod_{k=2}^n p_V(t_k - t_{k-1}, x_{k-1}, x_k) \, Z(T - t_n, x_n) \, d\mathbf{x} \text{,}
\end{equation*}
where
\begin{equation*}
    Z_V(t,x)
    = \int \limits_\R p_V(t,x,y) \, dy \text{,}
\end{equation*}
$0 < t_1 < \dots < t_n < T$, $B_1, \ldots, B_n \in \mathcal{B}(\R^n)$.

We use this formula to construct the penalizing measure in the case of the potential $\mu\, \delta(x - a)$. Let's start with the kernel of the operator $e^{t \A_\mu}$.
\begin{theorem}
\label{th61}
    The kernel $p_\mu(t,x,y)$ of the operator $e^{t\A_\mu}$, $t>0$, is given by
    \begin{align}
    \label{eq:semigroup_kernel}
    \begin{split}
        p_\mu(t,x,y)
        = p_0(t,x,y)
        &+ e^{\nu t} \, \frac{\psi_\nu(x-a) \psi_\nu(y - a)}{\| \psi_\nu \|_2^2} \\
        &+ \frac{1}{2 \pi i} \int \limits_{\gamma - i\infty}^{\gamma + i \infty} e^{\lambda t} \, \frac{\psi_\lambda(x - a) \psi_\lambda(y - a)}{(\lambda - \nu)(\psi_\nu, \psi_{\bar{\lambda}})} \, d\lambda \text{,}
    \end{split}
    \end{align}
    where $\gamma\in (0, \nu)$, and the function $p_0(t,x,y)$ is the kernel of the operator $e^{t \A}$,
    \begin{equation*}
        p_0(t,x,y) = \frac{1}{2 \pi} \int \limits_\R e^{-ip(x-y)} e^{-t \Psi(p)} \, dp \text{.}
    \end{equation*}
\end{theorem}
\begin{proof}
    The kernel $p_\mu(t,x,y)$ is connected to the kernel $r_\mu(\lambda,x,y)$ of $(\A_\mu - \lambda)^{-1}$ by the formula
    \begin{equation*}
        p_\mu(t,x,y)
        = - \frac{1}{2 \pi i} \int \limits_{\chi - i\infty}^{\chi + i \infty} e^{\lambda t} \, r_\mu(\lambda,x,y) \, d\lambda \text{,}
    \end{equation*}
    where $\chi$ is a constant satisfying $\chi > \nu$.
    
    Combining \eqref{eq:res_kernel} and this formula, we get
    \begin{align*}
        p_\mu(t,x,y)
        &= \frac{1}{2 \pi} \int \limits_\R e^{-ip(x-y)} e^{-t \Psi(p)} \, dp \\
        &+ \frac{1}{2 \pi i} \int \limits_{\chi - i\infty}^{\chi + i \infty} e^{\lambda t} \, \frac{\psi_\lambda(x - a) \psi_\lambda(y - a)}{(\lambda - \nu)(\psi_\nu, \psi_{\bar{\lambda}})} \, d\lambda \text{.}
    \end{align*}
     
     Shifting the contour $(\chi - i\infty, \chi +i\infty)$ to the left in the last integral, we get the sum of the integral over the contour $(\gamma - i\infty, \gamma +i\infty)$, $\gamma\in (0,\nu)$, and the residue at the point $\nu$, which gives us \eqref{eq:semigroup_kernel}.
\end{proof}

\begin{remark*}
    If $t = 0$, then the operator $e^{t\A_\mu}$ is the identity operator. Thus,
    \begin{equation*}
        p_\mu(0,x,y) = \delta(x - y) \text{.}
    \end{equation*}
\end{remark*}

Let us introduce the normalizing function $Z_\mu(t,x)$, $x \in \R$:
\begin{equation*}
    Z_\mu(t,x)
    = \int \limits_\R p_\mu(t,x,y) \, dy \text{.}
\end{equation*}
Using the Theorem \ref{th61} and the Lemma \ref{l24}, we obtain
\begin{equation*}
    Z_\mu(t,x)
    = 1
    + e^{\nu t} \, \frac{\psi_\nu(x-a)}{\nu \| \psi_\nu \|_2^2}
    + \frac{1}{2 \pi i} \int \limits_{\gamma - i\infty}^{\gamma + i \infty} \int \limits_\R e^{\lambda t} \, \frac{\psi_\lambda(x - a) \psi_\lambda(y - a)}{(\lambda - \nu)(\psi_\nu, \psi_{\bar{\lambda}})} \, d\lambda \, dy \text{.}
\end{equation*}

\begin{theorem}
\label{th62}
    Two following equalities hold.
    \begin{itemize}
        \Item[1.]
        \begin{equation*}
            \E e^{\mu L(t,x-a)}
            = \int \limits_\R p_\mu(t, x, y) \, dy
            =  \big( e^{t \A_\mu} \, 1 \big)(x) \text{.}
        \end{equation*}
        \Item[2.]
        \begin{equation*}
            \lim_{t \to \infty} e^{-\nu t} \, \E e^{\mu L(t,x-a)}
            = \frac{\psi_\nu(x-a)}{\nu \| \psi_\nu \|_2^2} \text{.}
        \end{equation*}
    \end{itemize}
\end{theorem}
\begin{proof}
    For $M > 0$, consider the function $\varkappa_M \in C_0(\R)$ such that
    \begin{equation*}
        \varkappa_M(x)
        =
        \begin{cases}
            1\text{,} & |x| < M \\
            0\text{,} & |x| > M+1 \text{.}
        \end{cases}
    \end{equation*}

    From Theorem \ref{th51} and continuity of the solution of the Cauchy problem \eqref{eq:cauchy_problem2}, it follows that
    \begin{equation*}
        \E \varkappa_M(\xi_x(t)) e^{\mu L(t,x-a)}
        = \int \limits_\R \varkappa_M(y) p_\mu(t,x,y) \, dy \text{.}
    \end{equation*}
    Using finiteness of
    \begin{equation*}
        \E e^{\mu L(t,x-a)} \quad\text{and} \quad \int \limits_\R p_\mu(t,x,y) \, dy
    \end{equation*}
    and setting $M \to \infty$, we get that 
    \begin{equation*}
        \E e^{\mu L(t,x-a)}
        = \int \limits_\R p_\mu(t,x,y) \, dy
        = \big( e^{t \A_\mu} \, 1 \big) (x) \text{.}
    \end{equation*}

     Therefore,
     \begin{equation*}
        e^{-\nu t} \, \E e^{\mu L(t,x - a)}
        = e^{-\nu t} Z_\mu(t,x)
        = \frac{\psi_\nu(x-a)}{\nu \| \psi_\nu \|_2^2}
        + Q(t,x) \text{,}
    \end{equation*}
    where $Q(t,x) \rightarrow 0$ as $t \to \infty$, and the proof is complete.
\end{proof}

Recall that the measure $\mathbf{P}_{T,x}$ of the process $\xi_x(t)$, $t \leq T$, is a measure on the sample paths $\Omega_{T,x} = \{ \omega \in \mathbb{D}([0,T], \R) \, | \, \omega(0) = x \}$, where $\mathbb{D}([0,T],\R)$ is a Skorokhod space, that is, a set of right-continuous with left limits real-valued functions on $[0,T]$ equipped with the Skorokhod distance.

Introduce the measure $\mathbf{Q}_{T,x}^\mu$, defined on cylindrical sets of $\Omega_{T,x}$ by the formula
\begin{equation}
    \mathbf{Q}_{T,x}^\mu \{ \omega(t_1) \in B_1, \dots, \omega(t_n) \in B_n \}
\end{equation}
\begin{equation*}
    = \frac{1}{Z_\mu(T,x)} \, \int \limits_{B_1} \dots \int \limits_{B_n} p_\mu(t_1, x, x_1) \prod_{k=2}^n p_V(t_k - t_{k-1}, x_{k-1}, x_k) \, Z_\mu(T - t_n, x_n) \, d\mathbf{x} \text{,}
\end{equation*}
where $0 < t_1 < \dots < t_n \leq T$, $B_1, \dots, B_n \in \mathcal{B}(\R)$.

\begin{lemma}
    The finite-dimensional distributions of $\mathbf{Q}_{T,x}^\mu$ are consistent, that is,
    \begin{align*}
        &\mathbf{Q}_{T,x}^\mu \{ \omega(t_1) \in B_1, \dots, \, \omega(t_{n-1}) \in B_{n-1}, \omega(t_n) \in \R \} \\
        &= \mathbf{Q}_{T,x}^\mu \{ \omega(t_1) \in B_1, \dots, \omega(t_{n-1}) \in B_{n-1} \} \text{,}
    \end{align*}
    where $0 < t_1 < \dots < t_n \leq T$, $B_1, \dots, B_{n-1} \in\mathcal{B}(\R)$.
\end{lemma}
\begin{proof}
    It is sufficient to show that
    \begin{equation*}
        \mathbf{Q}_{T,x}^\mu \{ \omega(t_1) \in B_1, \omega(t_2) \in \R \}
        = \mathbf{Q}_{T,x}^\mu \{ \omega(t_1) \in B_1 \} \text{.}
    \end{equation*}
    Using Theorem \ref{th62}, we obtain
    \begin{equation*}
        \mathbf{Q}_{T,x}^\mu \{ \omega(t_1) \in B_1, \omega(t_2) \in \R \}
    \end{equation*}
    \begin{equation*}
        = \frac{1}{Z_\mu(T,x)} \, \int \limits_{B_1} \int \limits_{\R} p_\mu(t_1, x, x_1) \, p_\mu(t_2 - t_1, x_1, x_2) \, Z_\mu(T - t_2, x_2) \, dx_1 \, dx_2
    \end{equation*}
    \begin{equation*}
        = \frac{1}{Z_\mu(T,x)} \, \int \limits_{B_1} p_\mu(t_1, x, x_1) \Bigg( \int \limits_{\R^2} p_\mu(t_2 - t_1, x_1, x_2) \, p_\mu(t_2 - t_1, x_1, y) \, dx_2 \, dy \Bigg) dx_1
    \end{equation*}
    \begin{equation*}
        = \frac{1}{Z_\mu(T,x)} \, \int \limits_{B_1} p_\mu(t_1, x, x_1) \Bigg( e^{(t_2-t_1) \A_\mu} \Bigg( \int \limits_{\R} p_\mu(T - t_2, \cdot, y) \, dy \Bigg) \Bigg)(x_1) \, dx_1
    \end{equation*}
    \begin{equation*}
        = \frac{1}{Z_\mu(T,x)} \, \int \limits_{B_1} p_\mu(t_1, x, x_1) \big( e^{(t_2-t_1) \A_\mu} \big( e^{(T-t_2) \A_\mu} \, 1 \big) \big)(x_1) \, dx_1
    \end{equation*}
    \begin{equation*}
        = \frac{1}{Z_\mu(T,x)} \, \int \limits_{B_1} p_\mu(t_1, x, x_1) \big( e^{(T-t_1) \A_\mu} \, 1 \big)(x_1) \, dx_1
    \end{equation*}
    \begin{equation*}
        = \frac{1}{Z_\mu(T,x)} \, \int \limits_{B_1} p_\mu(t_1, x, x_1) \, Z(T-t_1, x_1) \, dx_1
        = \mathbf{Q}_{T,x}^\mu \{ \omega(t_1) \in B_1 \} \text{,}
    \end{equation*}
    which completes the proof.
\end{proof}

By $\pi_\nu$ we denote the distribution 
\begin{equation*}
    \frac{\psi_\nu^2(x'-a)}{\| \psi_\nu \|_2^2} \, dx' \text{.}
\end{equation*}
    
For a family of distributions $\{ \mathbf{Q}_{T,x}^\mu \}_{T \geq 0}$, the following limit theorem holds.
\begin{theorem}
\label{th63}
    As $T \to \infty$, the densities of the finite-dimensional distributions of $\{ \mathbf{Q}_{T,x}^\mu \}$ converge pointwise and in $L_1(\R)$ to the densities of the corresponding finite-dimensional distributions of $\mathbf{P}_x^\mu$ -- the measure of a Markov process $\zeta(t)$, $t\geq 0$, with transition density
    \begin{equation*}
        \rho_\mu(t,x,y)
        = e^{-\nu t} \, \frac{p_\mu(t,x,y) \psi_\nu(y-a)}{\psi_\nu(x-a)}
    \end{equation*}
    and the invariant distribution of $\pi_\nu$.
\end{theorem}
\begin{proof}
    By the Theorem \ref{th51} and Theorem \ref{th52}, we have
    \begin{equation*}
        \int \limits_\R \rho_\mu(t,x,y) \, dy
        = \frac{e^{-\nu t}}{\psi_\nu(x-a)} \int \limits_\R \psi_\nu(y-a) \, p_\mu(t,x,y) \, dy
    \end{equation*}
    \begin{equation*}
        = \frac{e^{-\nu t}}{\psi_\nu(x-a)} \, \big( e^{t \A_\mu} \, \psi_\nu(\cdot - a) \big) (x)
        = \frac{e^{-\nu t}}{\psi_\nu(x-a)} \, \big( e^{\nu t} \, \psi_\nu(\cdot - a) \big) (x)
        = 1 \text{.}
    \end{equation*}
    
    Again using Theorem \ref{th52} and symmetricity of $p_\mu(t,x,y)$ with respect to spatial variables $x,y$, we obtain
    \begin{equation*}
        \int \limits_\R \rho_\mu(t,x,y) \, \pi_\nu(dx)
        = \int \limits_\R \rho_\mu(t,x,y) \, \frac{\psi_\nu^2(x-a)}{\| \psi_\nu \|_2^2} \, dx
    \end{equation*}
    \begin{equation*}
        = \int \limits_\R e^{-\nu t} \, \frac{\psi_\nu(x-a) \psi_\nu(y-a)}{\| \psi_\nu \|_2^2} \, p_\mu(t,x,y) \, dx
    \end{equation*}
    \begin{equation*}
        = \frac{\psi_\nu(y-a)}{\| \psi_\nu \|_2^2} \, \E \eta_\nu(t,x) \big|_{x=y} 
        = \frac{\psi_\nu^2(y-a)}{\| \psi_\nu \|_2^2} \text{.}
    \end{equation*}
    Therefore, $\rho_\mu(t,x,y)$ is indeed a transition density of some Markov process with the invariant distribution of $\pi_\nu$.

    Let us prove the pointwise convergence of densities. For $0 < t_1 < \dots < t_n < T$, $x_1, \dots, x_n \in \R$, we have
    \begin{equation*}
        \frac{1}{Z_\mu(T,x)} \, p_\mu(t_1, x, x_1) \dots p_\mu(t_n - t_{n-1}, x_{n-1}, x_n) \, Z_\mu(T - t_n, x_n)
    \end{equation*}
    \begin{equation*}
        =  \rho_\mu(t_1, x, x_1) \dots \rho_\mu(t_n - t_{n-1}, x_{n-1}, x_n) \, e^{\nu t_n} \, \frac{\psi_\nu(x-a) }{Z_\mu(T,x)} \, \frac{Z_\mu(T - t_n, x_n)}{\psi_\nu(x_n-a)} 
    \end{equation*}
    \begin{equation*}
        =  \rho_\mu(t_1, x, x_1) \dots \rho_\mu(t_n - t_{n-1}, x_{n-1}, x_n) \, \frac{1 + q_1(T,x,x_n)}{1 + q_2(T,x,x_n)} \text{,}
    \end{equation*}
    where $q_1(T,x,x_n) \rightarrow 0$, $q_2(T,x,x_n) \rightarrow 0$ for $T\to \infty$.

    Convergence in $L_1(\R)$ follows from the pointwise convergence by the Scheff\'e's lemma \cite[ch. 5, \S 5.10]{Williams}.
\end{proof}
\begin{corollary*}
    As $T \to \infty$, the finite-dimensional distributions of $\mathbf{Q}_{T,x}^\mu$ converge in total variation to the corresponding finite-dimensional distributions of $\mathbf{P}_x^\mu$.
\end{corollary*}

Let's describe the connection between the process $\zeta(t)$ and the semigroup $\{ \widetilde{U}_t \}$.
\begin{theorem}
\label{th64}
    The semigroup generated by the process $\zeta(t)$ coincides with $\{\widetilde{U}_t\}$.
\end{theorem}
\begin{proof}
    The kernel of the semigroup generated by $\zeta(t)$ is its transition density $\rho_\mu(t,x,y)$. Let $g \in C_0(\R)$. We have
    \begin{equation*}
        \int \limits_\R g(y) \rho_\mu(t,x,y) \, dy
        = \frac{e^{-\nu t}}{\psi_\nu(x - a)} \int \limits_\R g(y) \psi_\nu(y - a) p_\mu(t,x,y) \, dy
    \end{equation*}
    \begin{equation*}
        = \frac{e^{-\nu t}}{\psi_\nu(x - a)} \big( e^{t \A_\mu} [g \, \psi(\cdot - a)] \big)(x)
    \end{equation*}
    \begin{equation*}
        = \frac{e^{-\nu t} \, \E g(\xi_x(t)) \psi_\nu(\xi_x(t) - a) e^{\mu L(t,x - a)}}{\psi_\nu(x - a)} 
        = \frac{\E \eta_\nu(t,x) g(\xi_x(t))}{\psi_\nu(x - a)} = ( \widetilde{U}_t g )(x) \text{.}
    \end{equation*}
\end{proof}
\begin{corollary*}
    The process $\zeta(t)$ is a Feller process, and its sample paths belong to the Skorokhod space $\mathbb{D}([0,\infty),\R)$.
\end{corollary*}
\begin{proof}
    By Theorem \ref{th54}, the semigroup $\{\widetilde{U}_t\}$ is a Feller semigroup. Thus, so the process $\zeta(t)$ is a Feller process. As a consequence, the sample paths of $\zeta(t)$ belong to the space $\mathbb{D}([0,\infty),\R)$ \cite[ch. III, 2.6]{RogWil}.
\end{proof}

The following statement complements the Theorem \ref{th63} and its corollary.

\begin{theorem}
    As $T \to \infty$, the distributions $\{ \mathbf{Q}^\mu_{T,x} \}$ weakly converge to the probability distribution $\mathbf{P}_x^\mu$.
\end{theorem}
\begin{proof}
    Let $A \in \F_t$, $T > t$. We have
    \begin{equation*}
        \mathbf{Q}^\mu_{T,x}(A)
        = \frac{\E \mathbbm{1}_A(\omega) e^{\mu L(T,x)}}{\E e^{\mu L(T,x)}}
        = \frac{\E \, \E \big[ \mathbbm{1}_A(\omega) e^{\mu L(T,x)} | \F_t \big]}{\E e^{\mu L(T,x)}}
    \end{equation*}
    \begin{equation*}
        = \frac{\E \Big[ \mathbbm{1}_A(\omega) e^{\mu L(t,x)} \E e^{\mu L(T-s,y)} \, \big|_{y=\xi_x(t)} \Big]}{\E e^{\mu L(T,x)}}
    \end{equation*}
    \begin{equation*}
        = \frac{\E \Big[ \mathbbm{1}_A(\omega) e^{-\nu s} e^{\mu L(t,x)} e^{-\nu (T - s)} \E e^{\mu L(T-s,y)} \, \big|_{y=\xi_x(t)} \Big]}{e^{-\nu T} \E e^{\mu L(T,x)}} \text{.}
    \end{equation*}
    By the Theorem \ref{th62}, the last expression tends to
    \begin{equation*}
        \frac{\E \big[ \mathbbm{1}_A(\omega) e^{-\nu s} e^{\mu L(t,x)} \psi_\nu(\xi_x(t) - a) \big]}{\psi_\nu(x - a)}
        = \frac{\E \big[ \mathbbm{1}_A(\omega) \eta_\nu(t,x) \big]}{\psi_\nu(x - a)}
        = \mathbf{P}_x^\mu(A)
    \end{equation*}
    as $T \to \infty$.
\end{proof}

Eventually, we described the Feller process defined by the exponential attraction of the sample paths of the process $\xi(t)$ to the point $a$, and showed that this process is determined by the function $\psi_\nu(\cdot - a)$, which is the eigenfunction of the operator $\A_\mu$.

Let's prove one more limit theorem related to $\psi_\nu(\cdot - a)$. Consider the distribution $\mathbf{R}_{T,x}^\mu$ of the random variable $\omega(T)$. This is the distribution of the point to which an attracted sample path of the process $\xi_x(t)$ has come by the moment $T$.

\begin{theorem}
    As $T \to \infty$, the density $\mathbf{r}_{T,x}^\mu$ of the distribution $\mathbf{R}_{T,x}^\mu$ converges pointwise and in $L_1(\R)$ to the density of $\nu \psi_\nu(\cdot - a)$.
    
    Moreover,
    \begin{equation*}
        \| \mathbf{r}^\mu_{T,x} - \nu \psi_\nu(\cdot - a) \|_1
        \leq \frac{2(1 + z_\mu(T,x))}{Z_\mu(T,x)} \text{,}
    \end{equation*}
    where
    \begin{equation*}
        z_\mu(T,x) = \frac{1}{2 \pi i} \int \limits_{\gamma - i\infty}^{\gamma + i \infty} \int \limits_\R e^{\lambda T} \, \frac{\psi_\lambda(x - a) \psi_\lambda(x' - a)}{(\lambda - \nu)(\psi_\nu, \psi_{\bar{\lambda}})} \, d\lambda \, dx' \text{.}
    \end{equation*}
\end{theorem}
\begin{proof}
    We have
    \begin{equation*}
        \mathbf{r}_{T,x}^\mu(y)
        = \frac{p_\mu(T,x,y)}{Z_\mu(T,x)}
        = \frac{e^{\nu T} \big( \psi_\nu(x - a) \psi_\nu(y-a)/\| \psi_\nu \|_2^2 + q_1(T,x,y) \big)}{e^{\nu T} \big( \psi_\nu(x - a) / (\nu \| \psi_\nu \|_2^2) + q_2(T,x,y) \big)} \text{,}
    \end{equation*}
    where $q_1(T,x,y) \rightarrow 0$, $q_2(T,x,y)\rightarrow 0$ for $T\to \infty$. Therefore
    \begin{equation*}
        \lim_{T \to \infty} \mathbf{r}_{T,x}^\mu(y)
        = \nu \psi_\nu(y-a) \text{.}
    \end{equation*}
    
    Convergence in $L_1(\R)$ follows from the Scheff\'e's lemma \cite[ch. 5, \S 5.10]{Williams}.

    Furthermore,
    \begin{equation*}
        \| \mathbf{r}_{T,x}^\mu - \nu \psi_\nu(\cdot - a) \|_1
        = \int \limits_\R \bigg| \frac{p_\mu(T,x,y)}{Z_\mu(T,x)} - \nu \psi_\nu(y-a) \bigg| \, dy
    \end{equation*}
    \begin{equation*}
        = \frac{1}{Z_\mu(T,x)} \int \limits_\R \big| p_\mu(T,x,y) - \nu \psi_\nu(y-a) Z_\mu(T,x) \big| \, dy
    \end{equation*}
    \begin{equation*}
        \leq \frac{1}{Z_\mu(T,x)} \Bigg(
        \int \limits_\R p_0(T,x,y) \, dy
        + \frac{1}{2 \pi i} \int \limits_{\gamma - i\infty}^{\gamma + i \infty} \int \limits_\R e^{\lambda T} \, \frac{\psi_\lambda(x - a) \psi_\lambda(y - a)}{(\lambda - \nu)(\psi_\nu, \psi_{\bar{\lambda}})} \, d\lambda \, dy
    \end{equation*}
    \begin{equation*}
        + \int \limits_\R \nu \psi_\nu(y-a) \, dy
        + z_\mu(T,x) \int \limits_\R \nu \psi_\nu(y-a) \, dy \Bigg)
        \leq \frac{2(1 + z_\mu(T,x))}{Z_\mu(T,x)} \text{.}
    \end{equation*}
\end{proof}
\begin{corollary*}
    The distribution $\mathbf{R}_{T,x}^\mu$ converges in total variation to the distribution
    \begin{equation*}
        \nu \psi_\nu(y - a) \, dy \text{.}
    \end{equation*}
\end{corollary*}